\documentclass[11pt]{article}
\usepackage{CJK}
\usepackage{graphicx}
\usepackage{fancyhdr}
\usepackage{amsmath,amsthm,amssymb}
\usepackage{bm}
\usepackage{mathrsfs}
\usepackage{multirow}
\usepackage{indentfirst}
\usepackage{url}
\usepackage{stmaryrd}
\usepackage{color}
\usepackage{epstopdf}
\usepackage{natbib}
\usepackage[caption=false]{subfig}
\usepackage[utf8]{inputenc}
\usepackage[backref]{hyperref}
\usepackage[capitalize]{cleveref}
\newtheorem{theorem}{Theorem}[section]
\newtheorem{prop}{Proposition}[section]
\newtheorem{remark}{Remark}

\numberwithin{equation}{section}
\makeatletter
\newcommand{\Rmnum}[1]{\expandafter\@slowromancap\romannumeral #1@}
\makeatother
\newtheorem{corollary}{Corollary}[section]
\allowdisplaybreaks

\begin{document}

\title{\bf Asymptotic Distribution-Free Tests for Ultra-high Dimensional  Parametric Regressions via Projected Empirical Processes and $p$-value Combination
%Asymptotic Distribution Free Tests for Ultra-High-Dimensional Regression Models via Projected Empirical Processes and $p$-values Combination
\footnote{Corresponding author. Lixing Zhu's research was supported by the grants (NSFC12131006, NSFC12471276) from the National Natural Scientific Foundation of China and the grant (CI2023C063YLL) from the Scientific and Technological Innovation Project of China Academy of Chinese Medical Science.}} %The authors thank the editor, the associate editor and two referees for their constructive suggestions that led to the  improvement of an early manuscript.}}
\author{Falong Tan$^1$, Shan Tang$^2$, and Lixing Zhu$^{3*}$ \\~\\
{\small {\small {\it $^1$ Department of Statistics and Data Science, Hunan University, Changsha, China} }}\\
{\small {\small {\it $^2$ School of Mathematics and Statistics, Wuhan University of Technology, Wuhan, China} }}\\
{\small {\small {\it $^3$ Department of Statistics, Beijing Normal University at Zhuhai, Zhuhai, China}}}
}
\date{}
\maketitle

\begin{abstract}
This paper develops a novel methodology for testing the goodness-of-fit of sparse parametric regression models based on projected empirical processes and $p$-value combination, where the covariate dimension may substantially exceed the sample size. In such ultra-high dimensional settings, traditional empirical process-based tests often fail due to the curse of dimensionality or their reliance on the asymptotic linearity and normality of parameter estimators---properties that may not hold under ultra-high dimensional scenarios. To overcome these challenges, we first extend the classic martingale transformation to ultra-high dimensional settings under mild conditions and construct a Cram\'{e}r-von Mises type test based on a martingale-transformed, projected residual-marked empirical process for any projection on the unit sphere. The martingale transformation renders this projected test asymptotically distribution-free and enables us to derive its limiting distribution using only standard convergence rates of parameter estimators. While the projected test is consistent for almost all projections on the unit sphere under mild conditions, it may still suffer from power loss for specific projections. Therefore, we further employ powerful $p$-value combination procedures, such as the Cauchy combination, to aggregate $p$-values across multiple projections, thereby enhancing overall robustness. Furthermore, recognizing that empirical process-based tests excel at detecting low-frequency signals while local smoothing tests are generally superior for high-frequency alternatives, we propose a novel hybrid test that aggregates both approaches using Cauchy combination. The resulting hybrid test is powerful against both low-frequency and high-frequency alternatives. Detailed simulation studies and two real-data analyses are conducted to illustrate the effectiveness of our methodology in ultra-high dimensional settings.
%This paper proposes a two-step methodology for testing the goodness-of-fit of sparse parametric regression models where the covariate dimension may significantly exceed the sample size. In such ultra-high-dimensional regimes, traditional empirical process-based tests often fail due to the curse of dimensionality or the breakdown of asymptotic linearity and normality in parameter estimators. To overcome these challenges, we first extend the classic martingale transformation to high-dimensional settings under mild conditions. We construct a Cramér-von Mises type test based on a martingale-transformed, projected residual marked empirical process. This transformation renders this projected test asymptotically distribution-free and allows for valid inference using only standard convergence rates of the parameter estimators.While the projected test is consistent for almost all directions on the unit sphere, specific projections may suffer from power loss. To address this, our second step employs the Cauchy combination procedure to aggregate $p$-values across multiple projections, enhancing overall robustness. Furthermore, recognizing that empirical process-based tests excel at detecting low-frequency signals while local smoothing tests are superior for high-frequency alternatives, we propose a novel hybrid test that aggregates both approaches. Extensive simulations and two real-world data analyses demonstrate that our methodology maintains high power across diverse alternative scenarios and outperforms existing benchmarks.
\\

\end{abstract}

\newpage
%\baselineskip=16pt

%\newpage

\setcounter{equation}{0}
\section{Introduction}
This research is motivated by the problem of testing the goodness-of-fit of ultra-high dimensional regression models, where the dimension of covariates may substantially exceed the sample size. Consider the regression model:
\begin{equation}\label{1.1}
Y = m(X)+ \varepsilon,
\end{equation}
where $Y \in \mathbb{R}$ is the response, $X$ is the $p$-dimensional covariate vector, $m(\cdot)=E(Y|X=\cdot)$ is the unknown regression function, and $\varepsilon$ is the error term satisfying $E(\varepsilon|X)=0$. Our objective is to test whether the mean function $m(\cdot)$ belongs to some parametric class of functions $\mathcal{M}=\{ m(\cdot, \beta): \beta \in \Theta \subset \mathbb{R}^q \}$ in ultra-high dimensional settings.
%Therefore, the hypotheses can be formulized as
%\begin{eqnarray*}
%H_0: m(x) = m(x, \beta_0) \quad {\rm for \ some} \ \beta_0 \in \Theta \quad
%{\rm versus} \quad
%H_1: m(x) \neq m(x, \beta) \quad {\rm for \ all} \ \beta \in \Theta.
%\end{eqnarray*}

There is an extensive literature on goodness-of-fit testing for regression models in low dimensional settings when the dimension $p$ is considered to be fixed and smaller than the sample size $n$. One primary methodology for model checking, known as local smoothing tests, is based on nonparametric estimation of conditional moment restrictions $E[\varepsilon(\beta_0)|X]$ for some $\beta_0 \in \Theta$, where $\varepsilon(\beta_0) = Y-m(X, \beta_0)$. Examples include \cite{hardle1993}, \cite{zheng1996}, \cite{dette1999}, \cite{fan2001}, \cite{horowitz2001}, \cite{koul2004}, \cite{keilegom2008}, \cite{lavergne2008,lavergne2012}, \cite{guo2016}. These tests are usually asymptotically distribution-free and are particularly sensitive to high-frequency alternative models \citep{horowitz2001}. However, due to their reliance on nonparametric estimation, they usually suffer severely from the curse of dimensionality.
%Recently, Tan et al. (2025) proposed a new nonparametric smoothing-based test which can be applied to check the adequacy of high dimensional models, even when the covariate dimension substantially exceeds the sample size.
The other main type of test for model checking constructs test statistics based on empirical processes, which circumvents the nonparametric estimation of conditional moment restriction $E[\varepsilon(\beta_0)|X]$. See, for instance, \cite{Bierens1982}, \cite{stute1997}, \cite{Bierens1990}, \citet*{Stute1998a}, \cite{stute2002}, \cite{zhu2003}, \cite{escanciano2006b}, \cite{stute2008}, \cite{escanciano2018}, \cite{cuesta2019}, \cite{lu2020}, \cite{escanciano2024}.
These empirical process-based tests are usually powerful against low-frequency alternative models and can detect local alternatives at the parametric rate $n^{-1/2}$, which is the optimal detection rate in hypothesis testing.

However, empirical process-based tests typically require the asymptotic linearity or normality of parameter estimators to derive their limiting null distributions. To illustrate this, we consider the seminal paper by \cite{stute1997}. That paper showed that many goodness-of-fit tests for regressions are based on the empirical process $ \hat{S}_n(t) = n^{-1/2}\sum_{i=1}^n \varepsilon_i(\hat{\beta}) I(X_i \leq t)$, where $\varepsilon_i(\hat{\beta}) = Y_i-m(X_i, \hat{\beta})$ with $\hat{\beta}$ being a consistent estimator of $\beta$ and $\{(X_i, Y_i)\}_{i=1}^n $ is an i.i.d. sample with the same distribution as $(X, Y)$. Under some regularity conditions and fixed dimensional settings, \cite{stute1997} showed that under the null hypothesis,
\begin{equation}\label{stute-process}
\hat{S}_n(t) = S_n^0(t) + \sqrt{n}(\hat{\beta} - \beta_0)^{\top}M(t) + o_p(1),
\end{equation}
uniformly in $t$, where $M(t) = E[m'(X,\beta_0) I(X \leq t)]$ and $S_n^0(t)=n^{-1/2}\sum_{i=1}^n \varepsilon_i(\beta_0) I(X_i \leq t)$ with $\varepsilon_i(\beta_0) = Y_i-m(X_i, \beta_0)$. It is readily seen that the asymptotically linear expansion or normality of $\sqrt{n}(\hat{\beta} - \beta_0)$ are required to derive the limiting null distribution of $\hat{S}_n(t)$. However, this asymptotic property for the estimated parameter $\hat{\beta}$ in high dimensional settings, such as Lasso, post-Lasso, or their variants, may no longer hold. Consequently, empirical process-based tests that incorporate these estimation methodologies, without further transformation for the corresponding empirical processes, typically cannot be directly extended to ultra-high dimensional settings where the covariate dimension $p$ significantly surpasses the sample size $n$. Another critical issue is that most existing empirical process-based tests suffer from the curse of dimensionality because of data sparsity in high dimensional spaces; see \cite{escanciano2006b} and \cite{tan2025a} for more details on this issue. More recently, \cite{Shah2018} and \cite{Jankova2020} proposed two goodness-of-fit tests based on residual prediction and generalized residual prediction for high dimensional linear and generalized linear models, respectively, where the covariate dimension $p$ can be much larger than the sample size $n$. However, \cite{Shah2018} considered a goodness-of-fit test for regression models with fixed design, which is different from the present paper. \cite{Jankova2020} proposed a generalized residual prediction (GRP) test based on projected residuals $w^{\top}\hat{R}$, where $\hat{R}$ and $w$ are the estimated residual vector and projection on the unit sphere, respectively. Note that the GRP test only considers a specific aspect of the model misspecification; it may fail to capture all potential departures from the null hypothesis and loses power against certain alternatives in high dimensional settings.

The purpose of this paper is to develop a new goodness-of-fit test for regression models which can be applied in ultra-high dimensional scenarios and simultaneously mitigate the curse of dimensionality. Recall that empirical process-based tests typically require asymptotic linearity or normality of the parameter estimator $\hat{\beta}$ to derive their limiting distributions. Since this assumption for $\hat{\beta}$ in ultra-high dimensional settings may not hold, corresponding empirical process-based tests may not be applied to these settings either. Interestingly, we find that the classic martingale transformation \citep{Stute1998b} for model checking may be used to address this problem. This method can be traced back to \cite{khmaladze1981} for deriving a goodness-of-fit test of the cumulative distribution function; see also \cite{koul1999}, \cite{bai2001}, \cite{koenker2002}, \cite{bai2003},  \cite{khmaladze2004,khmaladze2009},  \cite{delgado2008}, \cite{tan2019a}, \cite{lu2020}, and \cite{tan2025a}, among many others. A martingale transformation, say $T$ for instance, is a linear operator such that the resulting test based on the martingale-transformed empirical process can be asymptotically distribution-free. More specifically, it eliminates the shift function, such as $M(t)$ in the decomposition of $\hat{S}_n(t)$, by setting $TM(t) \equiv 0$, and simultaneously ensures that the transformed process $TS_n^0(t)$ admits the same asymptotic properties as $S_n^0(t)$. This implies that $TS_n^0(t)$ and consequently $T\hat{S}_n(t)$ will be asymptotically distribution-free, with a limit of a Brownian motion in transformed time. Note that the martingale transformation eliminates the function $M(t)$, and thus the shift term $\sqrt{n}(\hat{\beta} - \beta_0)^{\top}M(t)$ in $\hat{S}_n(t)$ would also vanish. Consequently, the martingale-transformed empirical process would not involve the estimated parameter $\hat{\beta}$ and then can be applicable for ultra-high dimensional model checking. However, the classic martingale transformation introduced by \cite{Stute1998b} for model checking was designed for univariate covariates. \cite{khmaladze2004} and \cite{delgado2008} further extended the univariate martingale transformation to multivariate cases in fixed dimensional scenarios. Nevertheless, both of these methods involved multiple integrals with respect to $p$-dimensional covariates \citep{lu2020}, which makes them difficult to be extended to ultra-high dimensional scenarios. To our knowledge, whether the martingale transformation can be applied in the settings where the covariate dimension $p$ may substantially exceed the sample size remains an open problem in the literature.

In this paper, inspired by the dimension-reduction test proposed by \cite{stute2002}, we successfully extended the classic martingale transformation to ultra-high dimensional settings under mind conditions. Building on this extension, we propose a novel methodology for testing the goodness-of-fit of sparse regression models, where the covariate dimension $p$ may substantially exceed the sample size $n$. Our methodology includes two steps.

First, we use the martingale-transformed projected residual-marked empirical process to construct the test statistic for any given projection on the unit sphere in $\mathbb{R}^p$. We establish the limiting null distributions of the martingale-transformed process and its corresponding projected test statistic under mild conditions, even when the dimension $p$ grows exponentially with the sample size $n$. Under the alternative hypothesis, this projected test is consistent for almost all projections on the unit sphere with asymptotic power $1$. By projecting the high dimensional covariates $X$ to a one-dimensional space, the projected test can detect local alternatives departing from the null at the parametric rate $n^{-1/2}$ while significantly mitigating the curse of dimensionality. Moreover, the proposed martingale transformation involves only a univariate integral, making it easy to compute in practice, even when the dimension $p$ is much larger than the sample size $n$. Theoretically, since the martingale transformation eliminates the shift term arising from parameter estimation, the projected test not only becomes asymptotically distribution-free but also requires only the standard convergence rate, rather than the asymptotic linearity or normality, of parameter estimators to derive the asymptotic properties under both the null and alternative hypotheses.

Second, we employ standard combination methods to combine the projected test statistics from different projections to form our final test. Note that although the projected tests can be consistent for almost all projections on the unit sphere, they may still lose power for some unsuitably chosen projections. Therefore, to avoid possible power loss for certain projections, we adopt the Cauchy combination method \citep{liu2020} to aggregate the corresponding $p$-values of the projected tests from different projections to enhance power. It is worth mentioning that empirical process-based tests are typically more sensitive than local smoothing tests for low-frequency alternative models, while local smoothing tests are generally more powerful for high-frequency or oscillating models \citep{FanLi2000, horowitz2001}. In practice, since the underlying forms of regression models are typically unknown, it is desirable to have a testing procedure that can be powerful for both high-frequency and low-frequency alternative models. To this end, we further propose a novel hybrid test that aggregates the combined empirical process-based test with the local smoothing test proposed by \cite{tan2025b} via the Cauchy combination method. Since the Cauchy combination is primarily influenced by the smallest $p$-values, the hybrid test inherits the advantages of empirical process-based tests and local smoothing tests. Simulation studies also show that the hybrid test is powerful for both high-frequency and low-frequency alternatives.

The remainder of this paper is organized as follows. Section 2 develops the projected residual-marked empirical process and establishes its limiting null distribution. In Section 3, we extend the classic martingale transformation to ultra-high dimensional settings and construct the projected test based on the martingale-transformed  projected residual-marked empirical process. The limiting null distributions of the martingale transformation and the corresponding test statistic are also established in this Section. Section 4 presents the power analysis for the martingale transformation and the projected test statistic. In Section 5, we construct the combined projected tests and discuss the choice of projections for practical use. Section 6 presents simulation studies and two real data analyses to assess the finite sample performance of our tests. Section 7 contains concluding remarks and topics for future study. All proofs for the theoretical results are deferred to the Supplementary Material.

\section{Projected residual-marked empirical process}
In this paper, we focus on testing the goodness-of-fit of generalized linear models (GLMs), when the covariate dimension $p$ may significantly exceed the sample size $n$. Our method can also be extended to test the adequacy of more generalized models, such as quasi-GLMs, parametric multiple index models, etc. Recall that under the GLM settings,
we have $E(Y|X=x) = \mu(\beta_0^{\top}X)$ and $var(Y|X=x) = V(\mu(\beta_0^{\top}x))$ for some unknown parameter $\beta_0 = (\beta_0^{(1)}, \dots, \beta_0^{(p)})^{\top} \in \mathbb{R}^{p}$ and some inverse link function $\mu(\cdot)$. To illustrate our method, we restrict ourselves to testing the misspecification of the conditional mean function $m(x) = E(Y|X=x)$. Consequently, the null and alternative hypotheses become
\begin{eqnarray*}
&& H_0: \mathbb{P}\{m(X) = \mu(\beta_0^{\top}X) \} = 1, \quad {\rm for \ some } \ \beta_0 \in \Theta, \\
&& H_1: \mathbb{P}\{m(X) \neq \mu(\beta^{\top}X) \} >0 , \quad {\rm for \ any } \ \beta \in \Theta,
\end{eqnarray*}
where $\Theta$ is a compact set in $\mathbb{R}^{p}$. In high dimensional settings with $p \geq n$, as shown by \cite{Jankova2020}, if the design matrix ${\bf X} = (X_1, \dots, X_n)^{\top}$ is of full rank, there always exists a solution $\beta_0 \in \mathbb{R}^{p}$ of the system of linear equations $m(X_i) = \mu(\beta_0^{\top}X_i)$ for $i=1, \dots, n$. This implies that GLMs can never be misspecified in practice without any model structural assumption when $p \geq n$. A commonly used assumption in high dimensional scenarios is the sparsity of regression models. It is readily seen that the problem of model checking for sparse GLMs is reasonable when $p \geq n$. Therefore, we consider sparse regression models under both the null and alternative hypotheses throughout this paper.
%the conditional density of $Y$ given $X=x$ is $ f(y|x, \beta) = c(x, y) \exp\{ \frac{y\theta(x) - b(\theta(x))}{\sigma^2} \}$, where $c(x, y)$ is a positive function and $ \theta(x) = \beta^{\top}x $ with $\beta = (\beta_1, \dots, \beta_p)^{\top} \in \mathbb{R}^p$.

We introduce some notations that will be used below. For a vector $\beta \in \mathbb{R}^p$, let $\beta^{(j)}$ denote the $j$-th entry of $\beta$ and let $\|\beta\|_q = (\sum_{j=1}^p |\beta^{(j)}|^q)^{1/q} $ for $ q \in \mathbb{Z}^{+}$ and $\|\beta\|_0$ be the number of non-zero entries of $\beta$. Let $I \subset \{ 1, 2, \dots, p \}$ and let $\beta_I$ denote the vector containing only the entries of $\beta$ whose indices are in $I$. For a matrix $M \in \mathbb{R}^{n \times p}$, let $M_I$ be the matrix only with the columns of $M$ whose indices are in $I$ and let $M_{I^c}$ be the columns of $M$ with the indices in the complement of $I$. Let $S \subset \{1, \dots, p \}$ be the active set that contains the indices of the covariates $X= (X^{(1)}, \dots, X^{(p)})^{\top}$ truly related to the response $Y$. Under the null $H_0$, the true regression parameter $\beta_0$ is sparse, and the active set becomes $S=\{j: \beta_0^{(j)} \neq 0 \}$.

Our methodology for testing the goodness-of-fit of ultra-high dimensional regression models depends on the following result.

\begin{prop}\label{prop-prob}
{\rm (i)} Let $W \in \mathbb{R}$ and $X \in \mathbb{R}^p$ be random variables. It follows that
\begin{eqnarray*}
E[W|X]=0 \ a.s. &\Longleftrightarrow& E[W|\alpha^{\top}X] =0 \ a.s. \ \ {\rm for \ all} \  \alpha \in \mathcal{S}^{p-1},\\
E[W|\alpha^{\top}X]=0 \ a.s. &\Longleftrightarrow&  E[W I(\alpha^{\top}X \leq t)] \equiv 0 \ \ {\rm for \ all} \ t \in \mathbb{R},
\end{eqnarray*}
where $\mathcal{S}^{p-1} = \{\alpha \in \mathbb{R}^p: \|\alpha\|_2 = 1  \}$.

{\rm (ii)} Suppose that $E|W|^2 < \infty$, $E\| X \|_2^k < \infty$, and $ \sum_{k=1}^{\infty} (E\| X \|_2^k)^{-1/k} = \infty$. If we write
$ \mathcal{A} = \{ \alpha \in \mathbb{R}^{p}: E[W |\alpha^{\top}X] = 0 \ a.s. \} $,  then
\begin{eqnarray*}
\mathbb{P}\{ E[W|X]=0 \} = 1 \quad \Longleftrightarrow \quad \mathcal{A} \ {\rm has \ positive \ Lebesgue \ measure.}
\end{eqnarray*}
Moreover, if we write $\mathcal{A}_1=\{ \alpha \in \mathcal{S}^{p-1}: E[W|\alpha^{\top}X]=0 \ a.s.\}$, then
\begin{eqnarray*}
&& \mathbb{P}\{ E[W|X]=0 \} = 1 \quad  \Longleftrightarrow \quad \mathcal{L}(\mathcal{A}_1) = 1,  \\
&& \mathbb{P}\{ E[W|X] \neq 0 \} > 0 \quad  \Longleftrightarrow \quad \mathcal{L}(\mathcal{A}_1) = 0,
\end{eqnarray*}
where $\mathcal{L}$ denotes the uniform probability measure on the unit sphere $\mathcal{S}^{p-1} $.
\end{prop}

Proposition~\ref{prop-prob} (i) has been established in Lemma 2.1 of \cite{zhu1998}, Lemma 1 of \cite{escanciano2006a}, or Lemma 2.1 of \cite{lavergne2008}. \cite{patilea2016} and \cite{cuesta2019} derived similar results to the first part of Proposition~\ref{prop-prob}(ii) in the setting of functional data. The condition $\sum_{k=1}^{\infty} (E\| X \|_2^k)^{-1/k} = \infty $ in Proposition~\ref{prop-prob}(ii) is called the Carleman's condition, which can be satisfied if the random vector $X$ has a finite moment generating function around the neighborhood of zero; see \cite{cuesta2007} for more details. A detailed proof of Proposition~\ref{prop-prob} is provided in the Supplementary Material.

%Based on Lemma \ref{lemma-pro}, we can readily obtain the following result which will be crucial to our test statistic construction.
%
%\begin{corollary}\label{cor-prob}
%Suppose that the conditions in part {\rm (ii)} of Lemma \ref{lemma-pro} holds. If we write $\mathcal{S}^{p-1} = \{\alpha \in \mathbb{R}^p: \|\alpha\|_2 = 1 \}$ and $ \mathcal{A}_1= \{ \alpha \in \mathcal{S}^{p-1}: E[W|\alpha^{\top}X]=0 \ a.s. \} $, then
%\begin{eqnarray*}
%&& \mathbb{P}\{ E[W|X]=0 \} = 1 \quad  \Longleftrightarrow \quad \mathcal{L}(\mathcal{A}_1) = 1,  \\
%&& \mathbb{P}\{ E[W|X] \neq 0 \} > 0 \quad  \Longleftrightarrow \quad \mathcal{L}(\mathcal{A}_1) = 0,
%\end{eqnarray*}
%where $\mathcal{L}$ is the uniform probability measure on the unit sphere $\mathcal{S}^{p-1}$.
%\end{corollary}

We write $\varepsilon(\beta) = Y-\mu(\beta^{\top}X)$. Since the null hypothesis $H_0$ is tantamount to $\mathbb{P}\{ E[\varepsilon(\beta_0)|X]=0 \} =1 $ for some $\beta_0 \in \Theta$, it follows from Proposition \ref{prop-prob} that $H_0$ holds if and only if
$ \mathcal{L}\{ \alpha \in \mathcal{S}^{p-1}: E[\varepsilon(\beta_0)|\alpha^{\top}X]=0 \ a.s. \} = 1 $ for some $\beta_0 \in \Theta$.
Therefore, to test the null hypothesis $H_0$, we may first choose a projection $\alpha \in \mathcal{S}^{p-1}$ and then test the projected null hypothesis
$$ H_0^{\alpha}: \mathbb{P}\{ E[\varepsilon(\beta_0)|\alpha^{\top}X] = 0 \} =1, \quad {\rm for \ some} \ \beta_0 \in \Theta. $$
The principle behind this testing methodology is as follows. Under the null $H_0$, the projected null $H_0^{\alpha}$ also holds.
Under the alternative $H_1$, we have $\mathbb{P}\{ E[\varepsilon(\beta)|X] \neq 0 \} >0 $ for all $\beta \in \Theta$ and then
$ \mathcal{L}\{ \alpha \in \mathcal{S}^{p-1}: E[\varepsilon(\beta)|\alpha^{\top}X]=0 \ a.s. \} = 0 $. This implies that, under the alternative $H_1$, the projected null hypothesis $H_0^{\alpha}$ fails for $\mathcal{L}$-a.s. projections on $\mathcal{S}^{p-1}$.
Consequently, the null $H_0$ is $\mathcal{L}$-a.s. equivalent to
$ H_0^{\alpha}: \mathbb{P}\{ E[\varepsilon(\beta_0)|\alpha^{\top}X] = 0 \} =1 $ for some $\beta_0 \in \Theta$. We then construct the test statistics according to the projected null hypothesis $H_0^{\alpha}$.

For any given projection $\alpha \in \mathcal{S}^{p-1}$, it follows from Proposition \ref{prop-prob} (i) that $H_0^{\alpha}$ is equivalent to
\begin{eqnarray}\label{popu-process}
E[\varepsilon(\beta_0)I(\alpha^{\top}X \leq t)] \equiv 0 \ \ \forall t \in \mathbb{R}, \quad {\rm for \ some} \ \beta_0 \in \Theta.
\end{eqnarray}
Let $\{(X_i, Y_i)\}_{i=1}^n $ be an i.i.d. sample with the same distribution as $(X, Y)$ and let $\hat{\beta}$ be an estimator of $\beta_0$ under the GLM setting, such as a penalized estimator or its variants.
Motivated by (\ref{popu-process}), we propose a projected residual-marked empirical process as
\begin{eqnarray}\label{Rn_alpha}
\hat{R}_n^{\alpha}(t) = \frac{1}{\sqrt{n}} \sum_{i=1}^n \varepsilon_i(\hat{\beta}) I(\alpha^{\top}X_{i} \leq t),
\end{eqnarray}
where $\varepsilon_i(\hat{\beta}) = Y_i - \mu(\hat{\beta}^{\top}X_i)$.

The projected residual-marked empirical process $\hat{R}_{n}^{\alpha}(t)$ depends on the chosen projection $\alpha$. However, selecting appropriate projections is a non-trivial task. An unsuitably chosen projection can lead to a significant loss of power under alternative hypotheses, especially in ultra-high dimensional settings. Section 5 provides a detailed discussion on the selection of projections to ensure that our proposed tests can achieve good power.

\subsection{Limiting null distribution of $\hat{R}_n^{\alpha}(t)$}
To derive the asymptotic properties of $\hat{R}_n^{\alpha}(t)$ under $H_0$ in ultra-high dimensional settings, we introduce some notions and regularity conditions. A random variable $W \in \mathbb{R}$ is called sub-Weibull of order $\tau > 0$, if
$$ \|W\|_{\psi_{\tau}}:=\inf\{\eta>0: E \psi_{\tau}(|W|/\eta) \leq 1 \} < \infty,  $$
where $\psi_{\tau}(x) = \exp(x^{\tau}) - 1$ for $x \geq 0$.
Since $ \|W\|_{\psi_{\tau}} = \inf\{\eta>0: E\exp(|W|^{\tau}/\eta^{\tau}) \leq 2 \} $, it follows from the Markov inequality that if $W$ is sub-Weibull of order $\tau$, then
$ \mathbb{P}(|W| \geq t) \leq 2 \exp(-t^{\tau} / \|W\|_{\psi_{\tau}}^{\tau}) $ for all $ t \geq 0 $.
It is readily seen that sub-Gaussian and sub-exponential random variables are special cases of sub-Weibull distributions of $\tau=2$ and $\tau=1$, respectively. Note that the mean of a sub-Weibull random variable is not required to be zero. More detailed results on sub-Weibull distributions are elaborated in \cite{Vladimirova2020} and \cite{Kuchibhotla2022}.
We further define $\hat{S} = \{j: \hat{\beta}^{(j)} \neq 0\}$, $\hat{s} = |\hat{S}|$, and $s=|S|$. The notation $C$ in the following denotes a constant independent of $n$, which may be different for each appearance.

(A1) Under $H_0$, the estimator $\hat{\beta}$ satisfies $\|\hat{\beta} - \beta_0 \|_2 = O_p(\sqrt{\frac{s \log{p}}{n}})$ and $\| \hat{\beta} - \beta_0 \|_1 = O_p(\sqrt{\frac{s^2 \log{p}}{n}})$.

(A2) The random variables $\varepsilon(\beta_0)$, $\mu'(\beta_0^{\top}X)$, and $\mu''(\beta_0^{\top}X)$ are sub-Weibull of order $\tau =2$ with $\max\{
\| \varepsilon(\beta_0) \|_{\psi_2}, \| \mu'(\beta_0^{\top}X) \|_{\psi_2}, \|\mu''(\beta_0^{\top}X) \|_{\psi_2}\} \leq C < \infty $.
The covariates $X \in \mathbb{R}^p$ are centered and $\alpha^{\top}X$ is sub-Weibull of order $\tau =2$ with $ \| \alpha^{\top}X \|_{\psi_2} \leq C < \infty $ for all $\alpha \in \mathcal{S}^{p-1}$.

(A3) The link function $\mu(\cdot)$ admits third derivatives, and $\max\{ |\mu''(\beta^{\top}x)|, |\mu'''(\beta^{\top}x)| \} \leq F(x) $ for all $\beta \in \Theta \subset \mathbb{R}^p$ with $F(X)$ being sub-weibull of order $\tau \geq 1/3$ and $\|F(X)\|_{\psi_{\tau}} \leq C < \infty$.

Condition (A1) is satisfied by many commonly used penalized estimators under mild conditions, such as the GLM Lasso estimator, the GLM SCAD estimator, or their variants, when the underlying GLMs are correctly specified (the null hypothesis), see \cite{buhlmann2011} for instance.
The sub-weibull assumption of order $2$ (sub-Gaussian) in (A2) is typically imposed in the literature of high dimensional data analysis \citep{Wainwright2019}. It is used to bound the tail probability of the remainder process in the decomposition of $\hat{R}_n^{\alpha}(t)$, when the dimension $p$ may exceed the sample size $n$.
Condition (A3) is satisfied by many GLMs that are used in practice, such as Gaussian linear models, logistic regression models, and probit regression models. It is also used to control the convergence rate of the remainder of the process  $\hat{R}_n^{\alpha}(t)$.
The sub-Weibull order $\tau \geq 1/3$ of $F(X)$ is a technical condition which can be weakened if we impose a more restrictive condition on the divergence rate of $p$ in Theorem \ref{th-Rn} below.

\begin{theorem}\label{th-Rn}
Suppose that Conditions (A1)-(A3) hold. If $ s^2\log^3{(p \vee n)} = o(n)$, then under $H_0$, uniformly in $t \in \mathbb{R}$,
\begin{eqnarray}\label{Rn}
\hat{R}_n^{\alpha}(t)
&=&  \frac{1}{\sqrt{n}} \sum_{i=1}^n \varepsilon_i(\beta_0) I(\alpha^{\top}X_i \leq t) - \sqrt{n}
     (\hat{\beta} -\beta_0)^{\top}M(t) + o_p(1) \nonumber \\
&=:& R_{n1}^{\alpha}(t) + o_p(1),
\end{eqnarray}
where $p \vee n = \max\{p, n\}$, $\varepsilon_i(\beta_0) = Y_i - \mu(\beta_0^{\top}X_i)$, and $M(t) = E[X \mu'(\beta_0^{\top} X)I(\alpha^{\top} X \leq t)]$.
\end{theorem}

Note that the asymptotic linearity or normality may not hold for some penalized estimators $\hat{\beta}$, such as Lasso, post-Lasso, and their variants, making it challenging to establish the limiting null distribution of $\hat{R}^{\alpha}_n(t)$ for these estimators in high dimensional scenarios. Consequently, the conventional Cram\'{e}r-von Mises functional or the Kolmogorov-Smirnov functional of $\hat{R}^{\alpha}_n(t)$ cannot be directly applied to construct goodness-of-fit tests for regression models. To overcome this difficulty, we extend the classic method of martingale transformation to ultra-high dimensional settings under mild conditions. Since the martingale transformation would make the shift term $M(t)$ in (\ref{Rn}) to vanish, we can derive the limiting null distribution of the resulting martingale-transformed process, even in the absence of the asymptotic linear expansion of $\hat{\beta}_{0} - \beta_{0}$. It is also worth noting that \cite{cuesta2019} used projected empirical processes to test the goodness of fit of functional linear models. However, their test statistic is based on the projected residual-marked empirical process without martingale transformation, and consequently it still requires the asymptotic linear expansion of parameter estimators to derive its limiting distribution.

\section{Martingale transformation in ultra-high dimensional settings}
In this section, we construct the martingale transformation for the projected residual-marked empirical process $\hat{R}_n^{\alpha}(t)$, which can be applied in settings where the covariate dimension $p$ may significantly exceed the sample size $n$. It follows from Theorem \ref{th-Rn} that under $H_0$,
\begin{eqnarray*}
\hat{R}_n^{\alpha}(t) = R_{n0}^{\alpha}(t)  - \sqrt{n}(\hat{\beta} - \beta_0)^{\top} M(t) + o_p(1),
\end{eqnarray*}
where $ R_{n0}^{\alpha}(t) = n^{-1/2}\sum_{i=1}^n \varepsilon_i(\beta_0) I(\alpha^{\top}X_i \leq t)$ and $M(t) = E[X \mu'(\beta_0^{\top}X) I(\alpha^{\top}X \leq t)]$.
Straightforward calculations show that under $H_0$,
\begin{equation}
cov(R_{n0}^{\alpha}(s), R_{n0}^{\alpha}(t)) = \psi_n^{\alpha}(s \wedge t),
\end{equation}
where $\psi_n^{\alpha}(t) =E[\varepsilon^2(\beta_0) I(\alpha^{\top}X \leq t)] = \int_{-\infty}^t \sigma^2_{\alpha}(u) dF_{\alpha}(u)$ and $\sigma^2_{\alpha}(u)=E[\varepsilon^2(\beta_0)|\alpha^{\top}X = u]$.
Note that for any given projection $\alpha$, the function class $\{I(\alpha^{\top}x \leq t): t \in \mathbb{R} \}$ is a VC class with a VC-index 2. It is readily seen that the empirical process $R_{n0}^{\alpha}(t)$ is asymptotically tight and the convergence of the finite-dimensional distributions of $R_{n0}^{\alpha}(t)$ can be proved by standard arguments. This yields
\begin{equation*}
R_{n0}^{\alpha}(t) \longrightarrow B(\psi(t)), \quad {\rm in \ distribution}
\end{equation*}
in the Skorohod space $D[-\infty, \infty]$, where $\psi(t)$ is the pointwise limit of $\psi_n^{\alpha}(t)$ and $B(t)$ is the standard Brownian motion.

The purpose of martingale transformation is to eliminate $M(t)$ in the shift term of $\hat{R}_n^{\alpha}(t)$ and simultaneously transform $R^{\alpha}_{n0}(t)$ to an innovation process that admits the same limiting null distribution as $R_{n0}^{\alpha}(t)$.
The principle of martingale transformation is as follows.
Let $A^{\alpha}(t)=\frac{\partial M(t)}{\partial \psi_n^{\alpha}(t)}$ be the Radon-Nikodym derivative of $M(t)$ with respect to $\psi_n^{\alpha}(t)$. Recall that
$ M(t) = E[X \mu'(\beta_0^{\top} X)I(\alpha^{\top} X \leq t)] %= E[r_{\alpha}(\alpha^{\top}X) I(\alpha^{\top}X \leq t)]
= \int_{-\infty}^t r_{\alpha}(u) F_{\alpha}(du), $
where $r_{\alpha}(u) = E[X \mu'(\beta_0^{\top} X)|\alpha^{\top}X = u]$.
It follows that $A^{\alpha}(t) = r_{\alpha}(t)/\sigma^2_{\alpha}(t)$. We then define a matrix as
\begin{equation*}
\Gamma^{\alpha}(t)
= \int_t^{\infty} A^{\alpha}(u)A^{\alpha}(u)^{\top} d\psi_n^{\alpha}(u)
= \int_t^{\infty} \frac{r_{\alpha}(u) r_{\alpha}(u)^{\top}}{\sigma^4_{\alpha}(u)} d\psi_n^{\alpha}(u).
\end{equation*}
%It can also be stated as
%$\Gamma^{\alpha}(t) = E[r_{\alpha}(\alpha^{\top}X) X^{\top} \mu'(\beta_0^{\top}X) I(\alpha^{\top}X \geq t) /\sigma^2_{\alpha}(\alpha^{\top}X)]$.
Assuming $\Gamma^{\alpha}(t)$ is nonsingular for any $t \in \mathbb{R}$, the martingale transformation is defined as
\begin{equation}\label{martingale}
Tf(t) = f(t) -\int_{-\infty}^t A^{\alpha}(u)^{\top} \Gamma^{\alpha}(u)^{-1} \int_u^{\infty} A^{\alpha}(v) df(v) d\psi_n^{\alpha}(u),
\end{equation}
where $f(t)$ is either a bounded variation function or a stochastic process such that the integral in~(\ref{martingale}) is well defined. It may also be a vector of functions sometimes.
Note that $\Gamma^{\alpha}(t) = \int_t^{\infty} A^{\alpha}(u)dM(u)^{\top} $ and $ M(t)= \int_{-\infty}^t A(u) d\psi_n^{\alpha}(u)$, it is readily seen that
\begin{eqnarray*}
TM(t)^{\top}
= M(t)^{\top} - \int_{-\infty}^t A^{\alpha}(u)^{\top} \Gamma^{\alpha}(u)^{-1} \int_u^{\infty} A^{\alpha}(v) dM(v)^{\top}
    d\psi_n^{\alpha}(u)
%= M(t)^{\top} - \int_{-\infty}^t A^{\alpha}(u)^{\top} d\psi_n^{\alpha}(u)
\equiv 0.
\end{eqnarray*}
This implies that the martingale transformation $T$ eliminates the shift term $M(t)$ in the decomposition of $ \hat{R}_n(t)$.
Also note that $T$ is a linear operator, it follows that $TR_{n1}^{\alpha}(t) = TR_{n0}^{\alpha}(t)$.

We further investigate the asymptotic properties of $ TR_{n0}^{\alpha}(t)$ under the null hypothesis. Recalling that $R_{n0}^{\alpha}(t) = n^{-1/2}\sum_{i=1}^n \varepsilon_i(\beta_0) I(\alpha^{\top}X_i \leq t)$, it follows that
\begin{eqnarray*}
TR_{n0}^{\alpha}(t)
%&=& R_{n0}^{\alpha}(t) - \int_{-\infty}^t A^{\alpha}(u)^{\top} \Gamma^{\alpha}(u)^{-1} \int_u^{\infty}A^{\alpha}(v)dR_{n0}^{\alpha}(v) d\psi_n^{\alpha}(u)  \\
= R_{n0}^{\alpha}(t) - \frac{1}{\sqrt{n}} \sum_{i=1}^n \varepsilon_i(\beta_0) \int_{-\infty}^t A^{\alpha}(u)^{\top}
  \Gamma^{\alpha}(u)^{-1} A^{\alpha}(\alpha^{\top}X_i) I(\alpha^{\top}X_i \geq u) d\psi_n^{\alpha}(u).
\end{eqnarray*}
This implies that $TR_{n0}^{\alpha}(t)$ is a centered cusum process with the covariance function
\begin{equation}\label{cov-martin}
cov(TR_{n0}^{\alpha}(s), TR_{n0}^{\alpha}(t)) = cov(R_{n0}^{\alpha}(s), R_{n0}^{\alpha}(t)) = \psi_n^{\alpha}(s \wedge t).
\end{equation}
Therefore, the transformed process $TR_{n0}^{\alpha}(t)$ has the same covariance structure as $R_{n0}^{\alpha}(t)$ under $H_0$.
The assertion (\ref{cov-martin}) is justified in the Supplementary Material.
Similar to the arguments for proving Theorem 1.2 of \cite{Stute1998b}, we can derive the asymptotic tightness and the finite-dimensional convergence of $TR_{n0}^{\alpha}(t)$. Consequently,
$$ TR_{n1}^{\alpha}(t) = TR_{n0}^{\alpha}(t) \longrightarrow B(\psi(t)), \quad {\rm in \ distribution} $$
in the Skorohod space $D[-\infty, \infty)$. Furthermore, under mild conditions, we can show that
$ T\hat{R}^{\alpha}_n(t) - TR_{n1}^{\alpha}(t) = o_p(1) $
uniformly in $t$. Altogether we obtain that the martingale-transformed process $T\hat{R}^{\alpha}_n(t)$, after the time transformation $z= \psi(t)$, converges in distribution to the standard Brownian motion $B(t)$.

For practical application, the martingale transformation $T$ needs to be estimated by its empirical analog. For this, recall that
$$ A^{\alpha}(t) = \frac{r_{\alpha}(t)}{\sigma^2_{\alpha}(t)}  \quad  {\rm and}  \quad
\Gamma^{\alpha}(t) = E\left( \frac{r_{\alpha}(\alpha^{\top}X)r_{\alpha}(\alpha^{\top}X)^{\top}}{\sigma^2_{\alpha} (\alpha^{\top}X)} I(\alpha^{\top}X \geq t)\right), $$
where $r_{\alpha}(t) = E[X \mu'(\beta_0^{\top} X)|\alpha^{\top}X = t]$ and $\sigma^2_{\alpha}(t) =E[\varepsilon^2(\beta_0)|\alpha^{\top}X = t]$.
Since we do not make any assumption for the quantities $r_{\alpha}(t)$ and $\sigma^2_{\alpha}(t)$ other than smoothness, they should be estimated in a nonparametric way.
In low dimensional settings, standard nonparametric estimators such as the Nadaraya–Watson estimator may be suitable for this purpose.
However, we note that $r_{\alpha}(t) = E[X \mu'(\beta_0^{\top} X)|\alpha^{\top}X = t] \in \mathbb{R}^p$ is a $p$-dimensional function. When the dimension $p$ significantly exceeds the sample size $n$, nonparametric estimations for $r_{\alpha}(t)$ would bring enormous theoretical difficulties in deriving the asymptotic properties of the corresponding estimated martingale transformation. To address this difficulty, inspired by \cite{stute2002}, we proposed a novel variant of martingale transformation under mild conditions.

(A4) The covariates $X$ satisfy the linear conditional mean assumption, i.e., $E[X|B^{\top}X] = D B^{\top}X $ for a non-random matrix $D \in \mathbb{R}^{p \times 2}$, where $B=(\alpha, \beta_0) \in \mathbb{R}^{p \times 2}$. Furthermore, the $(i,j)$-entry $D_{ij}$ of $D$ satisfies $|D_{ij}| \leq C < \infty$ for all $1 \leq i \leq p$ and $j \in \{1,2\}$.

The linear conditional mean assumption (A4) is satisfied if the distribution of the covariates $X$ is elliptically symmetric, such as normal distributions. In high dimensional settings, as \cite{hall1993} demonstrated, if the original covariate dimension $p$ is large, then $E[X|B^{\top}X]$ is approximately linear in $X$.

Under Condition (A4), we have
\begin{eqnarray*}
M(t)
= E\{E[X|B^{\top}X] \mu'(\beta_0^{\top} X)I(\alpha^{\top}X \leq t)\}
= D E[B^{\top}X \mu'(\beta_0^{\top} X)I(\alpha^{\top}X \leq t)].
\end{eqnarray*}
We write $M_0(t) = E[B^{\top}X \mu'(\beta_0^{\top} X)I(\alpha^{\top}X \leq t)]$, it follows that
$$ M(t) = D M_0(t) = D \left( \int_{-\infty}^t u g_{1\alpha}(u) dF_{\alpha}(u), \int_{-\infty}^{t} g_{2\alpha}(u) dF_{\alpha}(u) \right)^{\top}, $$
where $g_{1\alpha}(u) = E[\mu'(\beta_0^{\top} X)|\alpha^{\top}X=u]$ and $g_{2\alpha}(u) = E[\beta_0^{\top}X \mu'(\beta_0^{\top} X)|\alpha^{\top}X=u]$.
Consequently, the residual marked empirical process $\hat{R}^{\alpha}_n(t)$ in (\ref{Rn}) can be restated as
$$ \hat{R}^{\alpha}_n(t) = R_{n0}^{\alpha}(t) - \sqrt{n} (\hat{\beta} - \beta_0)^{\top}D M_0(t) + o_p(1). $$
To eliminate the shift term $\sqrt{n} (\hat{\beta} - \beta_0)^{\top}D M_0(t)$, we respectively define the new $A^{\alpha}(t)$ and $\Gamma^{\alpha}(t)$ as
$$ A^{\alpha}(t)
= \frac{\partial M_0(t)}{\partial \psi_n^{\alpha}(t)}
= \left( \frac{tg_{1\alpha}(t)}{\sigma^2_{\alpha}(t)}, \frac{g_{2\alpha}(t)}{\sigma^2_{\alpha}(t)} \right)^{\top}
$$
and
\begin{eqnarray*}
\Gamma^{\alpha}(t)
= \int_t^{\infty} A^{\alpha}(u) A^{\alpha}(u)^{\top} d\psi_n^{\alpha}(u)
= \int_t^{\infty} \frac{1}{\sigma^4_{\alpha}(u)}
  \begin{pmatrix}
    u^2 g_{1\alpha}(u)^2, & ug_{1\alpha}(u)g_{2\alpha}(u) \\
    u g_{1\alpha}(u)g_{2\alpha}(u), & g_{2\alpha}(u)^2
  \end{pmatrix}
  d\psi_n^{\alpha}(u).
\end{eqnarray*}
It is important to note that these new $A^{\alpha}(t)$ and $\Gamma^{\alpha}(t)$ only involve univariate function $g_{1\alpha}(t)$, $g_{2 \alpha}(t)$, and $\sigma_{\alpha}^2(t)$. We then estimate these quantities by one-dimensional Nadaraya–Watson estimators:
\begin{eqnarray*}
\hat{\sigma}_{\alpha}^2(t) &=& \frac{\sum_{i=1}^n \varepsilon_i^2(\hat{\beta}) K_h(t-\alpha^{\top}X_i)}{\sum_{i=1}^n K_h(t-\alpha^{\top}X_i)}, \\
\hat{g}_{1\alpha}(t)      &=& \frac{\sum_{i=1}^n \mu'(\hat{\beta}^{\top}X_i) K_h(t-\alpha^{\top}X_i)}{\sum_{i=1}^n K_h(t-\alpha^{\top}X_i)}, \\
\hat{g}_{2\alpha}(t)      &=& \frac{\sum_{i=1}^n \hat{\beta}^{\top}X_i \mu'(\hat{\beta}^{\top}X_i)
                               K_h(t-\alpha^{\top}X_i)}{\sum_{i=1}^n K_h(t-\alpha^{\top}X_i)},
\end{eqnarray*}
where $\varepsilon_i(\hat{\beta}) = Y_i -\mu(\hat{\beta}^{\top}X_i)$ and $K_h(\cdot) = h^{-1} K(\cdot/h) $ with univariate kernel function $K(\cdot)$ and bandwidth $h$.
The estimators $\hat{A}^{\alpha}_n(t)$ and $\hat{\Gamma}^{\alpha}_n(t) $ for $A^{\alpha}(t)$ and $\Gamma^{\alpha}(t)$ respectively are
\begin{eqnarray}
\hat{A}^{\alpha}_n(t) = \left(\frac{t \hat{g}_{1\alpha}(t)}{\hat{\sigma}_{\alpha}^2(t)}, \frac{\hat{g}_{2\alpha}(t)}{\hat{\sigma}_{\alpha}^2(t)} \right)^{\top}
\quad {\rm and} \quad
\hat{\Gamma}^{\alpha}_n(t) = \int_t^{\infty} \hat{A}^{\alpha}_n(u) \hat{A}^{\alpha}_n(u)^{\top} d\hat{\psi}_n^{\alpha}(u),
\end{eqnarray}
where $\hat{\psi}_n^{\alpha}(u) = n^{-1} \sum_{i=1}^n \varepsilon_i^2(\hat{\beta}) I(\alpha^{\top}X_i \leq u)$.
%$\hat{F}_{\alpha}(\cdot)$ is the empirical distribution function of $\{\alpha^{\top}X_i \}_{i=1}^n$.
Consequently, we obtain the empirical analogue $\hat{T}_n$ of the martingale transformation $T$:
\begin{eqnarray}
\hat{T}_n \hat{R}^{\alpha}_n(t)
= \hat{R}^{\alpha}_n(t) - \int_{-\infty}^t \hat{A}^{\alpha}_n(u)^{\top} \hat{\Gamma}^{\alpha}_n(u)^{-1} \int_u^{\infty} \hat{A}^{\alpha}_n(v) d\hat{R}^{\alpha}_n(v) d\hat{\psi}_n^{\alpha}(u),
\end{eqnarray}
where $\hat{R}^{\alpha}_n(t) = n^{-1/2} \sum_{i=1}^n \varepsilon_i(\hat{\beta}) I(\alpha^{\top}X_{i} \leq t)$ is given by (\ref{Rn_alpha}).

\subsection{Limiting null distribution of $\hat{T}_n \hat{R}^{\alpha}_n(t)$}
To derive the asymptotic properties of the martingale-transformed process $\hat{T}_n \hat{R}^{\alpha}_n(t)$ in ultra-high dimensional settings, we impose some additional regularity conditions.

(A5) The kernel function $K(\cdot)$ satisfies (i) $K(\cdot)$ is continuous on $\mathbb{R}$ and has a continuous derivative on its support $[-1, 1]$;
(ii) $K(x)=K(-x)$ and $K(\cdot)$ is of bounded variation; (iii) $\int_{-1}^1 K(u)du = 1$ and $\int_{-1}^1 u^{i} K(u)du = 0$ for $i=1, \dots, k-1$.

(A6) We write $\sigma_{1\alpha}^2(t) = \sigma^2_{\alpha}(t)f_{\alpha}(t) $, $w_{1\alpha}(t) = g_{1\alpha}(t) f_{\alpha}(t)$, and $w_{2\alpha}(t) = g_{2\alpha}(t) f_{\alpha}(t)$. The functions $\sigma^2_{1\alpha}(t)$, $w_{1\alpha}(t)$, and $w_{2\alpha}(t)$ admit derivatives up to the $k-1$ order in $t$. Let $w_{1\alpha}^{(i)}(t) = \frac{d^{i} w_{1\alpha}(t)}{dt^i}$, $w_{2\alpha}^{(i)}(t) = \frac{d^{i} w_{2 \alpha}(t)}{dt^i}$, and  $(\sigma_{1\alpha}^2)^{(i)}(t) =\frac{d^{i}\sigma_{1\alpha}^2(t)}{dt^i}$ for $i=1, 2, \dots, k-1$.
The functions $w_{1\alpha}^{(k-1)}(t)$, $w_{2\alpha}^{(k-1)}(t)$, and $(\sigma_{1\alpha}^2)^{(k-1)}(t)$ satisfy the Lipschitz condition:
\begin{eqnarray*}
|w_{1\alpha}^{(k-1)}(t + u) - w_{1\alpha}^{(k-1)}(t)| &\leq& L|u|, \quad \forall \ u \in U,  \\
|w_{2\alpha}^{(k-1)}(t + u) - w_{2\alpha}^{(k-1)}(t)| &\leq& L|u|, \quad \forall \ u \in U,  \\
|(\sigma_{1\alpha}^2)^{(k-1)}(t + u) - (\sigma_{1\alpha}^2)^{(k-1)}(t)| &\leq& L|u|, \quad \forall \ u \in U,
\end{eqnarray*}
for some neighborhood $U$ of zero. Moreover, we assume that $\sup_{t \in \mathbb{R}}|f_{\alpha}(t)| \leq C < \infty$ and $ \inf_{t \in \mathbb{R}} \sigma^2_{1 \alpha}(t) \geq C > 0 $ for all $\alpha \in \mathcal{S}^{p-1}$.

(A7) The bandwidth $h$ satisfies $\sqrt{n}h^{2k} = o(1)$ and $\log^4{n} = o(nh^4)$ as $n \to \infty$.

(A8) The matrix $\Gamma^{\alpha}(t)$ satisfies $\inf_{t \leq t_0}|\det(\Gamma^{\alpha}(t))| > 0$ for any $t_0 \in \mathbb{R}$ and $\alpha \in \mathcal{S}^{p-1}$, where $\det(\Gamma^{\alpha}(t))$ denotes the determinant of the matrix $\Gamma^{\alpha}(t)$.

Conditions (A5)-(A7) are usually used in the literature of high-order nonparametric estimation; see, for instance, Chapters 2 and 4 of \cite{rao1983} and \cite{zhu1996}. Condition (A8) is necessary for the uniform boundedness of $ \| \Gamma^{\alpha}(t)^{-1} \|_2$ from infinity.

The next result establishes the asymptotic property of the martingale transformed process $\hat{T}_n \hat{R}^{\alpha}_n(t) $ under $H_0$. Its proof is provided in the Supplementary Material.

\begin{theorem}\label{TnVn}
Suppose that Conditions (A1)-(A8) hold and the matrix $\Gamma^{\alpha}(t)$ is nonsingular for all $t \in \mathbb{R}$. If $s^4 \log^4{(p \vee n)} = o(n)$ and $\log^5{p} = o(n)$ as $n \to \infty$, then under $H_0$,
$$ \hat{T}_n \hat{R}^{\alpha}_n(t) \longrightarrow B(\psi(t)), \quad in \ distribution $$
in the Skorohod space $D[-\infty, \infty)$.
\end{theorem}

It is worth mentioning that the martingale transformation can have a more appealing structure when the null hypothesis is a Gaussian linear model. We also require the linear conditional mean assumption of the covariates $X$.

(A4$^{\prime}$) The covariates $X$ satisfy the linear conditional mean assumption, i.e., $E[X|\alpha^{\top}X] = \alpha^{\top}X D_0 $ for some non-random vector $D_0 \in \mathbb{R}^p$ with $|D_0^{(j)}| \leq C < \infty$ for $1 \leq j \leq p$.

Recalling that $\mu'(t) \equiv 1$ under Gaussian linear models, it follows from (A4$^{\prime}$) and (\ref{Rn}) that
$$ M(t) = E[X I(\alpha^{\top}X \leq t)] = E\{ E[X|\alpha^{\top}X] I(\alpha^{\top} X \leq t) \} = E[\alpha^{\top}X I(\alpha^{\top}X \leq t)] D_0. $$
We write $M_0(t) = E[\alpha^{\top}X I(\alpha^{\top}X \leq t)] = \int_{-\infty}^t u dF_{\alpha}(u) $.
Since $\sigma^2_{\alpha}(t) \equiv \sigma^2$ is a constant under Gaussian linear models, it follows that the quantities $A^{\alpha}(t)$ and $\Gamma^{\alpha}(t)$ can be restated as
\begin{eqnarray*}
A^{\alpha}(t) = \frac{\partial M_0(t)}{\partial \psi_n^{\alpha}(t)} =\frac{t}{\sigma^2},
\ \ {\rm and}  \ \
\Gamma^{\alpha}(t) = \int_t^{\infty} A^{\alpha}(u)^2 d\psi_n^{\alpha}(u) = \frac{1}{\sigma^4} \int_t^{\infty} u^2 d\psi_n^{\alpha}(u),
\end{eqnarray*}
where $ \psi_n^{\alpha}(t) = \int_{-\infty}^t \sigma^2_{\alpha}(u) dF_{\alpha}(u) = \sigma^2 F_{\alpha}(t)$.
Consequently, the martingale transformation $T$ can also be applied in this much simplified scenario.
Note that under Gaussian linear models, we avoid the nonparametric estimations for $\sigma^2_{\alpha}(u)=E[\varepsilon^2(\beta_0)|\alpha^{\top}X = u]$, $g_{1\alpha}(u) = E[\mu'(\beta_0^{\top} X)|\alpha^{\top}X=u]$ and $g_{2\alpha}(u) = E[\beta_0^{\top}X \mu'(\beta_0^{\top} X)|\alpha^{\top}X=u]$ when constructing the estimated martingale transformation. Therefore, the resulting martingale transformation can have a much simpler structure when testing ultra-high dimensional Gaussian linear models.

\subsection{Projected test statistics based on $\hat{T}_n \hat{R}^{\alpha}_n(t)$}
According to Proposition~\ref{prop-prob} and Theorem~\ref{TnVn}, we can construct the test for $H_0$ based on any functional of $\hat{T}_n \hat{R}^{\alpha}_{n}(t)$ such that the resulting test is asymptotically distribution-free. Specifically, we employ the Cram\'{e}r-von Mises functional of the martingale-transformed process $\hat{T}_n \hat{R}^{\alpha}_n(t)$ to construct the test statistic. We define the (informal) test statistic as:
\begin{equation*}
CvM_{n,\alpha}^2 = \int_{-\infty}^{t_0} |\hat{T}_n \hat{R}^{\alpha}_{n}(t)|^2 d\hat{\psi}_n^{\alpha}(t),
\end{equation*}
where $\hat{\psi}_n^{\alpha}(t) = n^{-1} \sum_{i=1}^n \varepsilon_i^2(\hat{\beta}) I(\alpha^{\top}X_i \leq t)$. By Theorem~\ref{TnVn} and the Extended Continuous Mapping Theorem \citep[Theorem 1.11.1]{van1996}, we have under $H_0$,
$$ CvM_{n,\alpha}^2  \longrightarrow  \int_{-\infty}^{t_0}  B^2(\psi(t)) d\psi(t)   \quad  {\rm in \ distribution}, $$
where $B(t)$ is a standard Brownian motion.
Since $B(t \psi(t_0))/\sqrt{\psi(t_0)}=B(t)$ in distribution, it follows that
$$ \int_{-\infty}^{t_0}  B^2(\psi(t)) d\psi(t) = \psi^2(t_0) \int_{0}^{1} B^2(t) dt, \quad  {\rm in \ distribution}. $$
Therefore, our final test statistic based on a given projection $\alpha$ is
\begin{equation}\label{3.12}
TCvM_{n, \alpha}^2 = \frac{1}{\hat{\psi}_n^{\alpha}(t_0)^2} \int_{-\infty}^{t_0} |\hat{T}_n\hat{R}_n^{\alpha}(t)|^2 d\hat{\psi}_n^{\alpha}(t),
\end{equation}
where $\hat{\psi}_n^{\alpha}(t_0)=n^{-1} \sum_{i=1}^{n} [Y_i-\mu(\hat{\beta}^{\top} X_i)]^2 I(\alpha^{\top}X_i \leq t_0)$ is an empirical analog of $\psi(t_0)$.
Applying Theorem~\ref{TnVn} and the Extended Continuous Mapping Theorem again, we readily obtain the limiting null distribution of $TCvM_{n, \alpha}^2$ in ultra-high dimensional settings.

\begin{corollary}\label{null-MCVM}
Assume the conditions of Theorem~\ref{TnVn}. Then, under $H_0$,
\begin{equation}\label{limit-Mcvm}
 TCvM_{n, \alpha}^2  \longrightarrow  \int_{0}^{1} B^2(t)dt  \quad {\rm in \ distribution}.
\end{equation}
\end{corollary}

Corollary~\ref{null-MCVM} implies that the projected test $TCvM_{n, \alpha}^2$ based on the martingale-transformed process $\hat{T}_n\hat{R}^{\alpha}_{n}(t)$ is asymptotically distribution-free, and consequently its critical values can be tabulated.

In homoscedastic cases, such as Gaussian linear models, we have $\sigma_{\alpha}^2(t) \equiv \sigma^2$, which is independent of $t$. Recall that $\sigma_{\alpha}^2(t) = E[\varepsilon^2(\beta_0)|\alpha^{\top}X=t] $ and $\psi_n^{\alpha}(t)=E[\varepsilon^2(\beta_0)I(\alpha^{\top}X \leq t)]=\sigma^2 F_{\alpha}(t)$.
Therefore, $\psi_n^{\alpha}(t)$ can be estimated by $\hat{\sigma}^2 \hat{F}_{\alpha}(t)$, where $ \hat{\sigma}^2= n^{-1} \sum_{i=1}^{n} [Y_i-\mu(\hat{\beta}^{\top} X_i)]^2 $ and $\hat{F}_{\alpha}(t)$ is the empirical distribution of $\{\alpha^{\top}X_i: i=1, \dots, n\}$.
Consequently, the test statistic $TCvM_{n,\alpha}^2$ becomes
$$ TCvM_{n,\alpha}^2 = \frac{1}{\hat{\sigma}^2 \hat{F}_{\alpha}(t_0)^2} \int_{-\infty}^{t_0} |\hat{T}_n \hat{R}_{n}^{\alpha}(t)|^2 d\hat{F}_{\alpha}(t). $$
For $t_0$, we adopt the $99\%$ quantile of $\hat{F}_{\alpha}$ for practical applications, as suggested by \cite{Stute1998b} and \cite{stute2002}.

\section{Power analysis}
In this section, we investigate the asymptotic properties of the martingale-transformed process $\hat{T}_n \hat{R}^{\alpha}_n(t)$ and the projected test $TCvM_{n,\alpha}^2$ under various alternative hypotheses. Consider the following alternative hypotheses converging to the null hypothesis $H_0$ at the rate $\gamma_n = n^{-c}$:
$$ H_{1n}:  m(x) = E(Y|X=x) = \mu(\beta_0^{\top}x) + \gamma_n L(x), $$
where $c \in [0, 1/2]$, and $L(\cdot)$ is a non-constant function with $\mathbb{P}\{L(X) = 0\} < 1$.
Here, $c = 0$ corresponds to the global alternative hypothesis and $c > 0$ corresponds to local alternative hypotheses. We also assume the sparsity for the regression function $m(\cdot)$.

To derive the asymptotic properties of $\hat{T}_n \hat{R}^{\alpha}_n(t)$ under various alternative hypotheses in ultra-high dimensional settings, we impose some additional conditions. Let
$\tilde{\sigma}^2_{\alpha}(t) = E[\varepsilon^2(\tilde{\beta}_0)|\alpha^{\top}X=t], \tilde{g}_{1\alpha}(t) = E[\mu'(\tilde{\beta}_0^{\top}X)|\alpha^{\top}X=t]$,
$\tilde{g}_{2\alpha}(t) = E[\tilde{\beta}_0^{\top}X  \mu'(\tilde{\beta}_0^{\top}X)|\alpha^{\top}X=t] $, and $l_{2\alpha}(t) = E[L(X)^2|\alpha^{\top}X=t]$.

%Recall that $\hat{S} = \{j: \hat{\beta}_0^{(j)} \neq 0\}$, and $\tilde{S} = \{j: \tilde{\beta}_0^{(j)} \neq 0\}$ with $\tilde{\beta}_0 = (\tilde{\beta}_0^{(1)}, \dots, \tilde{\beta}_0^{(p)})^{\top}$ being the maximizer of $E[\rho(Y, X^{\top}\beta)]$ under the global alternative $H_1$.

(A9) (i) Under the global alternative $H_1$, there exists a parameter $\tilde{\beta}_0 = (\tilde{\beta}_0^{(1)}, \dots, \tilde{\beta}_0^{(p)})^{\top} \in \Theta$ such that
$$ \|\hat{\beta} - \tilde{\beta}_0\|_1 = O_p(\tilde{s} \sqrt{\frac{\log{p}}{n}}) \quad {\rm and} \quad
\|\hat{\beta}_0 - \tilde{\beta}_0\|_2 = O_p(\sqrt{\frac{\tilde{s} \log{p}}{n}}), $$
where $\tilde{s} = |\tilde{S}|$ and $\tilde{S} = \{j: \tilde{\beta}_0^{(j)} \neq 0 \}$.
(ii) Under the local alternative $H_{1n}$ with $c \in (0, 1/2]$, we have $\tilde{\beta}_{0} - \beta_{0} = r_n M^L + o_p(r_n)$ with $M^L_{\tilde{S}_1^c} = 0$ and $\| M^L \|_2 = O(1)$, and
$$ \|\hat{\beta}_0 - \beta_0\|_1 = O_p(\tilde{s} \sqrt{\frac{\log{p}}{n}} + \gamma_n \sqrt{\tilde{s}_1}) \quad {\rm and} \quad
\|\hat{\beta}_0 - \beta_0\|_2 = O_p(\tilde{s} \sqrt{\frac{\log{p}}{n}} + \gamma_n), $$
where $\tilde{s}_1 = |\tilde{S}_1|$, $\tilde{s} = |\tilde{S}|$, and $\tilde{S}_1 = S \cup \tilde{S}$ with $S = \{j: \beta_0^{(j)} \neq 0\}$.

(A10) The random variables $\varepsilon(\tilde{\beta}_0)$, $\mu'(\tilde{\beta}_0^{\top}X)$, and $\mu''(\tilde{\beta}_0^{\top}X)$ are sub-Weibull of order $\tau =2$ with $\max\{
\| \varepsilon(\tilde{\beta}_0) \|_{\psi_2}, \| \mu'(\tilde{\beta}_0^{\top}X) \|_{\psi_2}, \|\mu''(\tilde{\beta}_0^{\top}X) \|_{\psi_2}\} \leq C < \infty $. The random variable $L(X)$ is sub-Weibull of order $\tau \geq 2/3$ with $\| L(X) \|_{\psi_{\tau}} \leq C < \infty$.

(A11) We write $\tilde{\sigma}^2_{1\alpha}(t) = \tilde{\sigma}^2_{\alpha}(t) f_{\alpha}(t)$, $\tilde{w}_{1\alpha}(t) = \tilde{g}_{1\alpha}(t)f_{\alpha}(t)$, and $\tilde{w}_{2\alpha}(t) = \tilde{g}_{2\alpha}(t)f_{\alpha}(t)$. The functions $\tilde{\sigma}^2_{1\alpha}(t)$, $\tilde{w}_{1\alpha}(t)$, and $\tilde{w}_{2\alpha}(t)$ admit $k-1$ order derivative in $t$ and let $(\tilde{\sigma}_{1\alpha}^2)^{(i)}(t) =\frac{d^{i} \tilde{\sigma}_{1\alpha}^2(t)}{dt^i}$, $\tilde{w}_{1\alpha}^{(i)}(t) = \frac{d^{i} \tilde{w}_{1\alpha}(t)}{dt^i}$, and $\tilde{w}_{2\alpha}^{(i)}(t) = \frac{d^{i} \tilde{w}_{2 \alpha}(t)}{dt^i}$ for $i=1, 2, \dots, k-1$.
The functions $(\tilde{\sigma}_{1\alpha}^2)^{(k-1)}(t), w_{1\alpha}^{(k-1)}(t)$, $w_{2\alpha}^{(k-1)}(t)$, and $l_{2\alpha}(t)$ satisfy the Lipschitz condition
\begin{eqnarray*}
|w_{1\alpha}^{(k-1)}(t + u) - w_{1\alpha}^{(k-1)}(t)| &\leq& L|u|  \quad \forall \ u \in U,  \\
|w_{2\alpha}^{(k-1)}(t + u) - w_{2\alpha}^{(k-1)}(t)| &\leq& L|u|  \quad \forall \ u \in U,  \\
|(\tilde{\sigma}_{1\alpha}^2)^{(k-1)}(t + u) - (\tilde{\sigma}_{1\alpha}^2)^{(k-1)}(t)| &\leq& L|u|  \quad \forall \ u \in U, \\
|l_{2\alpha}(t+u)f_{\alpha}(t+u) - l_{2\alpha}(t)f_{\alpha}(t)| &\leq& L|u|  \quad \forall \ u \in U,
\end{eqnarray*}
for some neighborhood $U$ of zero. Moreover, we assume that $ \inf_{t \in \mathbb{R}} \tilde{\sigma}^2_{1 \alpha}(t) \geq C > 0 $ for all $\alpha \in \mathcal{S}^{p-1}$.

(A12) Let $\tilde{\Gamma}^{\alpha}(t) = \int_t^{\infty} \tilde{A}^{\alpha}(u) \tilde{A}^{\alpha}(u)^{\top} d\tilde{\psi}_n^{\alpha}(u) $, where $ \tilde{\psi}_n^{\alpha}(v) = E[\varepsilon(\tilde{\beta}_0)^2I(\alpha^{\top}X \leq v)] $ and $\tilde{A}^{\alpha}(t) = \left(\frac{t \tilde{g}_{1\alpha}(t)} {\tilde{\sigma}_{\alpha}^2(t)}, \frac{\tilde{g}_{2\alpha}(t)} {\tilde{\sigma}_{\alpha}^2(t)}\right)^{\top} $.
The matrix $\tilde{\Gamma}^{\alpha}(t)$ satisfies $\inf_{t \leq t_0}|\det( \tilde{\Gamma}^{\alpha}(t) )| \geq C >0$ for any $t_0 \in \mathbb{R}$, where $\det( \tilde{\Gamma}^{\alpha}(t) )$ denotes the determinant of the matrix $\tilde{\Gamma}^{\alpha}(t)$.

We show in the Supplementary Material that the GLM lasso estimator satisfies Condition (A9) under both the local and global alternative hypotheses (the misspecified models).
It is also worth noting that \cite{Buhlmann2015} showed that under misspecified linear models, the support $\tilde{S} $ of $\tilde{\beta}_0$ satisfies $\tilde{S} \subset S$ if the covariates $X$ follow a Gaussian distribution with positive definite covariance matrix.
Under the misspecified generalized linear model with fixed dimension $p$, \cite{lu2012} proved that $\tilde{S} = S$ if the true underlying model is also a generalized linear model with a misspecified link function and the linear conditional mean assumption is satisfied for $X$, that is, $E(\beta^{\top}X|\beta_0^{\top}X ) = b \beta_0^{\top}X + a$ for all $\beta \in \mathbb{R}^p$ and $a, b \in \mathbb{R}$. Both of these results provide evidence for Condition (A9).
Conditions (A10) and (A11) are similar to (A2) and (A3) in Section 3, which are used to control the convergence rate of the remainders of $\hat{T}_n \hat{R}_n^{\alpha}(t)$ under the alternative hypotheses. Condition (A12) is needed to ensure the uniform boundedness of $ \| \tilde{\Gamma}^{\alpha}(t)^{-1} \|_2$ away from infinity.

The next theorem establishes the asymptotic properties of $\hat{T}_n \hat{R}_n^{\alpha}(t)$ under various alternative hypotheses. Its proof is provided in the Supplementary Material.
We write $\tilde{G}^{\alpha}(t) = E[\varepsilon(\tilde{\beta}_0) I(\alpha^{\top}X \leq t)]$, $S_L^{\alpha}(t) = E[L(X)I(\alpha^{\top}X \leq t)]$, $\tilde{T}\tilde{G}^{\alpha}(t) = \tilde{G}^{\alpha}(t) -\int_{-\infty}^t \tilde{A}^{\alpha}(u)^{\top} \tilde{\Gamma}^{\alpha}(u)^{-1} \int_u^{\infty} \tilde{A}^{\alpha}(v) d\tilde{G}^{\alpha}(v) d\tilde{\psi}_n^{\alpha}(u) $, and
$ TS_L^{\alpha}(t) = S_L^{\alpha}(t) -\int_{-\infty}^t A^{\alpha}(u)^{\top} \Gamma^{\alpha}(u)^{-1} \int_u^{\infty} A^{\alpha}(v)dS_L^{\alpha}(v) d\psi_n^{\alpha}(u) $.

\begin{theorem}\label{TnVn-alternative}
Suppose that Conditions (A2)-(A12) holds. \\
{\rm (1)} Under $H_1$, if $\tilde{\Gamma}^{\alpha}(t)$ is non-singular for all $t \in \mathbb{R}$, and $\tilde{s}^2 \log^4{(p \vee n)} = o(n)$ as $n \to \infty$, then
$$ n^{-1/2} \hat{T}_n \hat{R}_n^{\alpha}(t) \longrightarrow L_1(t), \quad in \ probability, $$
where $L_1(t)$ is the pointwise limit of $\tilde{T}\tilde{G}^{\alpha}(t)$.  \\
{\rm (2)} Under $H_{1n}$ with $r_n = n^{-a}$ and $a \in (0, 1/2)$, if $\Gamma^{\alpha}(t)$ is non-singular for all $t \in \mathbb{R}$, $\gamma_n \tilde{s}_1^2 \log^2{(p \vee n)} = o(1)$, $\log^5{p} = o(n)$, and $\sqrt{n}\gamma_n^2 = o(h)$, then
$$ (n r_n^2)^{-1/2} \hat{T}_n \hat{R}_n^{\alpha}(t) \longrightarrow L_2(t), \quad in \ probability, $$
where $L_2(t)$ is  the pointwise limit of $TS_L^{\alpha}(t)$. \\
{\rm (3)} Under $H_{1n}$ with $r_n = n^{-1/2}$, if $\Gamma^{\alpha}(t)$ is non-singular for all $t \in \mathbb{R}$, $\tilde{s}_1^2 \log^2{(p \vee n)} = o(\sqrt{n})$, and $\log^5{p} = o(n)$, then
$$ \hat{T}_n \hat{R}_n^{\alpha}(t) \longrightarrow B(\psi(t)) + L_2(t), \quad in \ distribution, $$
in the Skorohod space $D[-\infty, \infty)$, where $B(\psi(t))$ is given in Theorem~\ref{TnVn}.
\end{theorem}

The following asymptotic result for the projected test statistic $TCvM_{n, \alpha}^2$ is a consequence of Theorem~\ref{TnVn-alternative} and the Extended Continuous Mapping Theorem \citep[Theorem 1.11.1]{van1996}.
We write $\tilde{\psi}_n^{\alpha}(t) = E[\varepsilon^2(\tilde{\beta}_0)I(\alpha^{\top}X \leq t)] =
E[\tilde{\sigma}^2(\alpha^{\top}X) I(\alpha^{\top}X \leq t)] $.

\begin{corollary}\label{TCVM-alternative}
Suppose that Conditions (A2)-(A12) holds. \\
{\rm (1)} Under $H_1$, if $\tilde{\Gamma}^{\alpha}(t)$ is non-singular for all $t \in \mathbb{R}$, and $\tilde{s}^2 \log^4{(p \vee n)} = o(n)$ as $n \to \infty$, then
$$ \frac{1}{n} TCvM_{n, \alpha}^2 \longrightarrow \frac{1}{\tilde{\psi}(t_0)} \int_{-\infty}^{t_0} |L_1(t)|^2 d\tilde{\psi}(t) \quad in \ probability, $$
where $ \tilde{\psi}(t)$ and $L_1(t)$ are the pointwise limits of $\tilde{\psi}_n^{\alpha}(t)$ and $\tilde{T}\tilde{G}^{\alpha}(t)$, respectively.  \\
{\rm (2)} Under $H_{1n}$ with $r_n = n^{-c}$ and $c \in (0, 1/2)$, if $\Gamma^{\alpha}(t)$ is non-singular for all $t \in \mathbb{R}$, $\gamma_n \tilde{s}_1^2 \log^2{(p \vee n)} = o(1)$, $\log^5{p} = o(n)$, and $\sqrt{n}\gamma_n^2 = o(h)$, then
$$\frac{1}{n r_n^2} TCvM_{n, \alpha}^2 \longrightarrow \frac{1}{\psi(t_0)} \int_{-\infty}^{t_0} |L_2(t)|^2 d\psi(t) \quad in \ probability, $$
where $\psi(t)$ and $L_2(t)$ are the pointwise limits of $\psi_n^{\alpha}$ and $TS_L^{\alpha}(t)$, respectively. \\
{\rm (3)} Under $H_{1n}$ with $r_n = n^{-1/2}$, if $\Gamma^{\alpha}(t)$ is non-singular for all $t \in \mathbb{R}$, $\tilde{s}_1^2 \log^2{(p \vee n)} = o(\sqrt{n})$, and $\log^5{p} = o(n)$, then
$$ TCvM_{n, \alpha}^2 \longrightarrow \int_0^1 |B(t) + L_2(\psi^{-1}(t \psi(t_0)))/\sqrt{\psi(t_0)}|^2 dt, \quad in \ distribution, $$
where $B(t)$ is a standard Brownian motion..
\end{corollary}

Corollary~\ref{TCVM-alternative} implies that if $L_1(t)$ and $L_2(t)$ are non-zero functions, then the projected test $TCvM_{n,\alpha}^2$ is consistent under the global alternative hypothesis and can detect the local alternatives distinct from the null at the parametric rate $n^{-1/2}$, even when the covariate dimension $p$ grows exponentially with the sample size $n$.

\section{The test statistics for practical use}
\subsection{Combined projected test statistics}
Note that a sequence of test statistics can be constructed based on various projections. According to Proposition~\ref{prop-prob} and Corollary~\ref{TCVM-alternative}, if the limit function $L_1(t)$ is non-zero, then the projected test $TCvM_{n,\alpha}^2$ is consistent under $H_1$ for almost all projections $\alpha \in \mathcal{S}^{p-1}$ with respect to the uniform measure $\mathcal{L}$. However, the test may still suffer from a substantial power loss for some inappropriate choice of projection, particularly in ultra-high dimensional settings. Another potential problem is that the value of the projected test statistic $TCvM_{n,\alpha}^2$ may vary across distinct projections, potentially leading to unstable power performance when relying on a single projection. To address these limitations, we propose to combine various projected test statistics $TCvM_{n,\alpha}^2$ to form a final test statistic, thereby enhancing the overall statistical power.

A variety of methods are available in the literature for combining test statistics or their corresponding $p$-values, such as the classic Fisher’s combination method \citep{fisher1925} and the more recent Cauchy combination method \citep{liu2020}.
For a set of chosen projections $\{ \alpha_i \in \mathcal{S}^{p-1}, i =1, \dots, k_1\}$, the asymptotic $p$-value of the projected test $TCvM_{n,\alpha_i}^2$ is given by
$\hat{p}_{1\alpha_i} = 1-\Psi(TCvM_{n,\alpha_i}^2) $,
where $\Psi(\cdot)$ denotes the cumulative distribution function of the random variable $\int_{0}^{1} B^2(t)dt$.
Similar to the arguments for Theorem \ref{TnVn} and by applying the Continuous Mapping Theorem, we have, under $H_0$,
\begin{equation}\label{limit-pvalue}
 (\hat{p}_{1\alpha_1}, \dots, \hat{p}_{1\alpha_{k_1}}) \longrightarrow (p_1, \dots, p_{k_1}), \quad {\rm in \ distribution},
\end{equation}
where $\{p_i: i=1, \dots, k_1\}$ are uniform random variables on $(0,1)$. Note that these $p$-values $p_1, \dots, p_{k_1}$ may be mutually dependent. Consequently, the classic Fisher's combination method cannot be applied here as it requires the independence between the $p$-values. In contrast, the Cauchy combination method is robust to dependence of the $p$-values; see \cite{liu2020} for more details on this issue. We therefore employ it to combine the $p$-values to form our final test statistic. The Cauchy combination-based test statistic is given by
\begin{equation}\label{cauchy}
TCvM^2_C = \sum_{i=1}^{k_1} w_i \tan\{ (\frac{1}{2} - \hat{p}_{1\alpha_i}) \pi \} = \sum_{i=1}^{k_1} w_i \tan\{ [\Psi(TCvM_{n,\alpha_i}^2)-\frac{1}{2}] \pi \},
\end{equation}
where $w_i$ are non-negative weights satisfying $\sum_{i=1}^{k_1} w_i = 1$. In this paper, we simply use the equal weights, i.e., $w_i=1/k_1$ for $i=1, \dots, k_1$, which performs very well in our simulation studies.
By (\ref{limit-pvalue}) and applying the Continuous Mapping Theorem again, we have, under $H_0$,
\begin{equation*}
TCvM^2_C \longrightarrow \sum_{i=1}^{k_1} w_i \tan\{ (\frac{1}{2} - p_i) \pi \}, \quad {\rm in \ distribution}.
\end{equation*}
Let $P_w = \sum_{i=1}^{k_1} w_i \tan\{ (\frac{1}{2} - p_i) \pi \}$.
If $p_1, \dots, p_{k_1}$ are i.i.d. uniform random variables on $(0,1)$, then it is readily seen that $P_w$ follows a standard Cauchy distribution, $Cauchy(0,1)$.
Furthermore, even when the $p$-values $p_1, \dots, p_{k_1}$ are mutually dependent, \cite{liu2020} demonstrated that $ \lim_{t \to \infty} \mathbb{P}(P_w > t) /\mathbb{P}(P_0 > t) =1$ under the null hypothesis, where $P_0$ denotes a standard Cauchy random variable. This implies that even when there exist dependencies between the $p$-values, the tail probability of $P_w$ can be approximated by that of the standard Cauchy distribution. Consequently, the asymptotic critical values of the Cauchy combination-based test $TCvM^2_C$ can also be approximated by quantiles of the standard Cauchy distribution.

In model checking, it is well known that empirical process-based tests, such as $TCvM_{n,\alpha}^2$, are more sensitive to low-frequency alternative models, while local smoothing tests are usually more sensitive and powerful against high-frequency alternative models. In practice, however, researchers usually do not know in advance which kind of models the underlying regression model comes from if no prior information is available. Therefore, it is important to propose a test statistic that can be sensitive to both low-frequency and high-frequency alternatives. Note that the Cauchy combination is most determined by the smallest $p$-values \citep{liu2020}; thus, a natural idea is to employ this combination method to combine empirical process-based tests and local smoothing tests to achieve this goal.
Recently, \cite{tan2025b} proposed a new local smoothing test for ultra-high dimensional regression models via projections, which also performs very well against high-frequency alternative models.
Their test statistic based on a given projection $\alpha$ is
\begin{equation}\label{PLS-alpha}
PLS_{n,\alpha}
= \frac{\sum_{1\leq i \neq j \leq n} \varepsilon_i(\hat{\beta}) \varepsilon_j(\hat{\beta})
  K(\frac{\alpha^{\top}X_i - \alpha^{\top}X_j}{h})}{\left(2 \sum_{1\leq i \neq j \leq n} \varepsilon_i^2(\hat{\beta}) \varepsilon_j^2(\hat{\beta}) K^2(\frac{\alpha^{\top}X_i - \alpha^{\top}X_j}{h}) \right)^{1/2}},
\end{equation}
where $K(\cdot)$ is a univariate kernel function and $h$ is the bandwidth.
According to Corollary 3.1 of \cite{tan2025b}, for any given projection $\alpha$, $PLS_{n, \alpha}$ converges in distribution to a standard normal distribution under $H_0$. Let $\hat{p}_{2\alpha_i} = 1-\Phi(PLS_{n,\alpha_i})$ be the asymptotic $p$-value of $PLS_{n,\alpha_i}$ for $i=1, \dots, k_2$, where $\Phi(\cdot)$ denotes the cumulative distribution function of the standard normal distribution. We then propose a hybrid test statistic based on the Cauchy combination method as
\begin{equation}\label{hybrid-test}
Hybrid_C = \sum_{i=1}^{k_1} w_i \tan\{ (\frac{1}{2} - \hat{p}_{1\alpha_i}) \pi \} + \sum_{j=1}^{k_2} v_j \tan\{ (\frac{1}{2} - \hat{p}_{2\alpha_j})\pi \},
\end{equation}
where the weights $w_i$ and $v_j$ satisfy $\sum_{i=1}^{k_1} w_i + \sum_{j=1}^{k_2} v_j = 1$. The asymptotic critical values of $Hybrid_C$ can also be determined by the quantiles of the standard Cauchy distribution.

\begin{remark}
The combined tests $TCvM^2_C$ and $Hybrid_C$ are both asymptotically distribution-free, and thus we do not need to resort to the resampling methods such as the wild bootstrap to approximate the limiting null distribution. Therefore, these tests are easy to implement in practice, particularly in ultra-high dimensional settings. Under the global alternative $H_1$ and mild conditions, it follows from Corollary~\ref{TCVM-alternative} and Theorem 3.3 of \cite{tan2025b} that the asymptotic $p$-values $\hat{p}_{1\alpha_i}$ and $\hat{p}_{2\alpha_i}$ converge to zero for almost all projections $\alpha_i \in \mathcal{S}^{p-1}$. Consequently, the combined test statistics $TCvM^2_C$ and $Hybrid_C$ diverge to infinity for almost all projections as the sample size $n \to \infty$. \cite{liu2020} theoretically demonstrated that the Cauchy combination of $p$-values is robust to dependent $p$-values and is most influenced by the smallest $p$-values. Therefore, even if the test statistics $TCvM_{n,\alpha_i}^2$ or $PLS_{n, \alpha_i}$ based on a single projection $\alpha_i$ may not be consistent, the combined test statistics $TCvM^2_C$ and $Hybrid_C$ can still exhibit robust power performance under the alternative hypothesis. We also note that our empirical process-based test $TCvM_{n, \alpha}^2$ is more sensitive to low-frequency alternatives and the local smoothing test $PLS_{n,\alpha}$ is more powerful against high-frequency models. Since the Cauchy combination of $p$-values is primarily determined by the smallest $p$-values, the hybrid test $Hybrid_C$ is expected to be powerful for both low-frequency and high-frequency alternative models. Simulation studies in Section 6 validate these theoretical assertions.
\end{remark}

\subsection{The choice of projections}
In practice, the power of the combined tests $TCvM^2_C$ and $Hybrid_C$ would be heavily influenced by the choice of projections in ultra-high dimensional scenarios. To illustrate this, we consider the problem of testing the goodness-of-fit of Gaussian linear models. Assume without loss of generality that the data are standardized, so we have $E(X) = 0$ and $E(Y) = 0$.
According to Corollary~\ref{TCVM-alternative}, the projected test $TCvM_{n, \alpha}^2$ can have power under the alternative hypothesis for a given projection $\alpha \in \mathcal{S}^{p-1}$ only if $E[\varepsilon(\tilde{\beta}_0) I(\alpha^{\top}X \leq t)] \neq 0$ for some $t \in \mathbb{R}$, where $\varepsilon(\tilde{\beta}_0) = Y-\tilde{\beta}_0^{\top}X$. Consider an extreme case where $\varepsilon(\tilde{\beta}_0)$ is mean independent of $\alpha^{\top}X$ for a certain projection $\alpha \in \mathcal{S}^{p-1}$, i.e., $E[\varepsilon(\tilde{\beta}_0)|\alpha^{\top}X] = E[\varepsilon(\tilde{\beta}_0)]$. This can be achieved if $X \sim N(0, I_p)$, $\varepsilon $ is independent of $X$, and $\alpha_S =0$, where $S$ is the true active set. Indeed, if $X \sim N(0, I_p)$, it follows from \cite{Buhlmann2015} that the support $\tilde{S} $ of $\tilde{\beta}_0$ satisfies $\tilde{S} \subset S$. Note that $\varepsilon(\tilde{\beta}_0) = m(X) + \varepsilon -\tilde{\beta}_0^{\top}X = m(X_S) + \varepsilon -\tilde{\beta}_{0S}^{\top} X_S $ and $\alpha^{\top}X = \alpha^{\top}_{S^c} X_{S^c} $, where $S^c$ denotes the complement of $S$. This implies that $\varepsilon(\tilde{\beta}_0)$ is independent of $\alpha^{\top}X$ and consequently $E[\varepsilon(\tilde{\beta}_0)|\alpha^{\top}X] = E[\varepsilon(\tilde{\beta}_0)]$.
Recalling that $E(X) = 0$ and $E(Y) = 0$, it follows that $E[\varepsilon(\tilde{\beta}_0)|\alpha^{\top}X] = 0$, and thus $E[\varepsilon(\tilde{\beta}_0) I(\alpha^{\top}X \leq t)] \equiv 0$. This means that the projected test $TCvM_{n, \alpha}^2$ may have no power under the alternative hypothesis in such scenarios. Therefore, to ensure good power, we should choose projections such that $\alpha^{\top}X$ has high correlation with the error term $\varepsilon(\tilde{\beta}_0)$. To this end, we propose a data-driven procedure to select the projections to ensure the proposed tests $TCvM^2_C$ and $Hybrid_C$ have good power under the alternatives.

Recall that $\varepsilon(\tilde{\beta}_0) = m(X) - \mu(\tilde{\beta}_0^{\top}X) + \varepsilon $, and we assume sparsity for the regression models under both the null and alternative hypotheses. Motivated by this, we assume without loss of generality that the error term $\varepsilon(\tilde{\beta}_0)$ admits a multiple-index model structure, that is, $\varepsilon(\tilde{\beta}_0) = g(\tilde{\vartheta}, \tilde{\theta}_1^{\top}X, \dots, \tilde{\theta}_d^{\top}X, \varepsilon) $, where $\tilde{\vartheta} \in \mathbb{R}$, and the projections $\tilde{\theta}_1, \dots, \tilde{\theta}_d \in \mathcal{S}^{p-1}$ are latent parameters. If there were no dimension reduction structure in this model, then we would have $d = p$ and $\tilde{\theta}_i = e_i$, where $e_i \in \mathcal{S}^{p-1}$ with $1$ in the $i$-th component and $0$ otherwise. However, this case is unlikely to occur given the sparsity assumed under both the null and alternative hypotheses. It is evident that the projected variables $\tilde{\theta}_1^{\top}X, \cdots, \tilde{\theta}_d^{\top}X$ are highly correlated with the error $\varepsilon(\tilde{\beta}_0)$. Furthermore, the projected predictor $\tilde{\theta}_0^{\top}X$ with $\tilde{\theta}_0 = \tilde{\beta}_0/\|\tilde{\beta}_0\|$ may also exhibit high correlation with the error \citep{stute2002}. Consequently, a natural idea is to construct the Cauchy combination-based test statistics based on these latent projections $\tilde{\theta}_0, \tilde{\theta}_1, \dots, \tilde{\theta}_d$.

However, all these latent projections are unknown and must be estimated in practice. If the full sample were used to estimate the latent projections, dependencies would be introduced between the estimated projections and the test statistics, thereby complicating the derivation of the asymptotic distribution. To address this issue, we randomly split the data into two parts, $\mathcal{D}_1$ and $\mathcal{D}_2$, of equal size. Here, we assume without loss of generality that the sample size $n$ is even. We use the first part of the data, $\mathcal{D}_1$, to estimate the projections $\tilde{\theta}_0, \tilde{\theta}_1, \dots, \tilde{\theta}_d$, and subsequently construct the test statistic based on the second part, $\mathcal{D}_2$. Let $\hat{\theta}_0^{(1)}, \hat{\theta}_1^{(1)}, \dots, \hat{\theta}_{\hat{d}^{(1)}}^{(1)}$ be the estimators of $\tilde{\theta}_0, \tilde{\theta}_1, \dots, \tilde{\theta}_d$ using $\mathcal{D}_1$, respectively. The resulting Cauchy combination test statistic is defined as:
\begin{equation}\label{cauchy_half}
TCvM^2_C
=\sum_{i=0}^{\hat{d}^{(1)}} w_i\tan\{ (\frac{1}{2} - \hat{p}^{(2)}_{1\hat{\theta}_i^{(1)}}) \pi \},
\end{equation}
where $\sum_{i=0}^{\hat{d}^{(1)}} w_i =1$, and $\hat{p}^{(2)}_{1\hat{\theta}_i^{(1)}} = 1 - \Psi(TCvM_{n,\hat{\theta}_i^{(1)}}^{2(2)})$ is computed using $\mathcal{D}_2$ and the estimated projection $ \hat{\theta}_i^{(1)} $.

The test $TCvM^2_C $ in~(\ref{cauchy_half}) uses only half of the data to construct the test statistic, which may lead to power loss under the alternative hypothesis. To address this problem, we adopt the cross-fitting strategy to enhance power. Specifically, we swap the roles of the two subsamples, $\mathcal{D}_2$ is used to estimate the latent projections, while $\mathcal{D}_1$ is employed to construct the test statistic. The resulting cross-fitting test statistic is given by
\begin{eqnarray}\label{cauchy_all}
TCvM^2_{CF}
&=& \sum_{i=0}^{\hat{d}^{(1)}} w_i\tan\{ (\frac{1}{2} -
    \hat{p}^{(2)}_{1\hat{\theta}_i^{(1)}}) \pi \}
    + \sum_{j=0}^{\hat{d}^{(2)}} w_j\tan\{ (\frac{1}{2} - \hat{p}^{(1)}_{1\hat{\theta}_j^{(2)}}) \pi \},
\end{eqnarray}
where $\sum_{i=0}^{\hat{d}^{(1)}} w_i + \sum_{j=0}^{\hat{d}^{(2)}} w_j = 1$, $\{ \hat{\theta}_j^{(2)} \}_{j=0}^{\hat{d}^{(2)}}$ are estimated using $\mathcal{D}_2$, and $\hat{p}^{(1)}_{1\hat{\theta}_j^{(2)}} = 1 - \Psi(TCvM_{n,\hat{\theta}_j^{(2)}}^{2(1)})$ is computed using $\mathcal{D}_1$ and the projection $\hat{\theta}_j^{(2)} $.
Similarly, the resulting cross-fitting hybrid test statistic is given by
\begin{eqnarray}\label{hybrid-dual}
Hybrid_{CF}
&=& \sum_{i=0}^{\hat{d}^{(1)}} w_i\tan\{ (\frac{1}{2} - \hat{p}^{(2)}_{1\hat{\theta}_i^{(1)}})
    \pi \} +\sum_{j=0}^{\hat{d}^{(2)}} w_j\tan\{ (\frac{1}{2} - \hat{p}^{(1)}_{1\hat{\theta}_j^{(2)}}) \pi \} \nonumber \\
&&  +\sum_{i=1}^{\hat{d}^{(1)}} v_i \tan\{ (\frac{1}{2} - \hat{p}^{(2)}_{2
    \hat{\theta}_i^{(1)}})\pi \}
    + \sum_{j=1}^{\hat{d}^{(2)}} v_j \tan\{ (\frac{1}{2} - \hat{p}^{(1)}_{2
    \hat{\theta}_j^{(2)}}) \pi\},
\end{eqnarray}
where the weights $w_i, w_j, v_i, v_j$ satisfy $\sum_{i=0}^{\hat{d}^{(1)}} w_i + \sum_{j=0}^{\hat{d}^{(2)}} w_j + \sum_{i=1}^{\hat{l}^{(1)}} v_i + \sum_{j=1}^{\hat{l}^{(2)}} v_j =1$, $\hat{p}^{(2)}_{2 \hat{\theta}_i^{(1)}} = 1- \Phi(PLS^{(2)}_{n, \hat{\theta}_i^{(1)}})$ is computed using $\mathcal{D}_2$ and the projection $\hat{\theta}_i^{(1)}$ from $\mathcal{D}_1$, and $\hat{p}^{(1)}_{2 \hat{\theta}_j^{(2)}} = 1- \Phi(PLS^{(1)}_{n, \hat{\theta}_j^{(2)}})$ is computed using $\mathcal{D}_1$ and the projection $\hat{\theta}_i^{(2)}$ from $\mathcal{D}_2$.

\section{Numerical studies}
\subsection{Simulations}
In this section, we conduct simulation studies to assess the finite sample performance of the proposed tests $TCvM^2_C $, $TCvM^2_{CF}$, and $Hybrid_{CF} $ when the covariate dimension $p$ may substantially exceed the sample size $n$. Since our tests rely on the parameter estimation, an inaccurate estimate of the parameters may affect the finite sample performance of the tests. We therefore employ the post-Lasso estimator of $\beta_0$ to construct the test statistics. The post-Lasso method applies least squares to the model selected by the Lasso estimator. According to \cite{Chernozhukov2013}, the post-Lasso estimator performs at least as well as Lasso in terms of the rate of convergence and has the additional advantage of a smaller bias. This implies that the post-Lasso estimator $\hat{\beta}_0$ satisfies Condition (A1) and the resulting tests are asymptotically distribution-free. We also conduct simulation studies for our tests using the Lasso estimator. These unreported results show that although our proposed tests control the empirical size in the setting of testing Gaussian linear models, they fail to maintain the significant level for testing logistic regression models. This may suggest that, in high dimensional scenarios, the finite sample performance of goodness-of-fit tests for regression models is also affected by the parameter estimation methodology.

To compute the tests $TCvM^2_C $, $TCvM^2_{CF}$, and $Hybrid_{CF}$, we randomly split the data into two parts, $\mathcal{D}_1$ and $\mathcal{D}_2$, and then construct the test statistics according to (\ref{cauchy_half}), (\ref{cauchy_all}), and (\ref{hybrid-dual}). To accurately estimate the latent projections $\tilde{\theta}_0, \tilde{\theta}_1, \dots, \tilde{\theta}_{d}$, we utilize the Distance Correlation Sure Independence Screening \citep{lirunze2012} to select the first $\frac{n/2}{\log(n/2)}$ top-ranked variables, and then apply the sparse sufficient dimension reduction technique LassoSIR \citep{lin2019} to construct the estimated projections $\hat{\theta}_i^{(1)}$ and $\hat{\theta}_i^{(2)}$ for $i=0,1, \dots, d$. Furthermore, since the proposed tests involve nonparametric estimation, we compute them using the Epanechnikov kernel $K(x) = (3/4)(1-x^2)I(|x| \leq 1)$ and choose the bandwidth by cross validation automatically.

We compare our tests with the $RP_n$ test of \cite{Shah2018}, the $GRP_n$ test of \cite{Jankova2020}, and the local smoothing test $\hat{T}_{\rm{Fisher}}^{C}$ proposed by \cite{tan2025b}. The latter is based on combining the projected tests $PLS_{n,\alpha}$ from different projections. For the specific forms of these tests, we refer the reader to the respective papers for details. In the following simulations, the parameter $a=0$ corresponds to the null hypothesis, and $a \neq 0$ corresponds to the alternative hypothesis. The simulation results are based on the average of 1000 replications with a nominal level of $\tau=0.05$. The simulation results of the $RP_n$ test are computed using the R package RPtests, and those for $GRP_n$ are computed by running the code available on the website \url{https://github.com/jankova/GRPtests}, provided by \cite{Jankova2020}. We compute the test $\hat{T}_{\rm{Fisher}}^{C}$ using the Epanechnikov kernel $K(x) = (3/4)(1-x^2)I(|x| \leq 1)$ and the bandwidth $h = n^{-2/9}$ as suggested in \cite{tan2025b}. Furthermore, despite all the following null models having zero intercepts, we estimate them as unknown parameters in the simulation studies.

In the first simulation study, we consider the case of testing Gaussian linear models, where the covariate dimension $p$ may be much larger than the sample size $n$.

{\em Study 1.} Generate data from the following models:
\begin{eqnarray*}
H_{11}:  Y &=& \beta_0^{\top}X+ 0.1a(\beta_0^{\top}X)^2 + \varepsilon;  \\
H_{12}:  Y &=& \beta_0^{\top}X+ a \cos(0.6 \pi \beta_0^{\top}X) + \varepsilon;\\
H_{13}:  Y &=& \beta_0^{\top}X+ a \exp(0.25 \beta_1^{\top}X) + \varepsilon;
\end{eqnarray*}
where $\beta_0 = (1,1,1,1,1,0,\dots,0)^{\top}$ and $\beta_1=(\underbrace{1,\dots,1}_{p_1},0,\dots,0)$ with $p_1 = 10$. The covariates $X \in \mathbb{R}^p$ is $N(0, \Sigma)$ independent of the standard Gaussian error $\varepsilon$, where $\Sigma=I_p$ or $\Sigma=(\rho^{|i-j|})_{p \times p}$ with $\rho=0.4$ and $\rho = 0.8$. Here, $H_{11}$ and $H_{13}$ are low-frequency models and $H_{12}$ is a high-frequency model under the alternative hypothesis. We consider the sample size $n=300$ with the covariate dimension $p \in \{50, 100, 300, 600, 900, 1200\}$.

The simulation results are presented in Tables 1-3. We observe that all tests control the significant level very well in all cases of sample sizes and dimensions. The $GRP_n$ test of \cite{Jankova2020} is slightly conservative with smaller empirical sizes. For low-frequency models $H_{11}$ and $H_{13}$, our tests $TCvM^{2}_{C}$, $TCvM^{2}_{CF}$, and $Hybrid_{CF}$ generally exhibit higher power than the other three competing tests across all settings. For the high-frequency model $H_{12}$, the tests $RP_n$, $GRP_n$, $TCvM^{2}_{C}$, and $TCvM^{2}_{CF}$ exhibit almost no power, even in the low dimensional setting ($p=50$). In contrast, the local smoothing test $\hat{T}_{\rm{Fisher}}^C$ and our test $Hybrid_{CF}$ have substantially higher power than the other four. This result is inline with the traditional findings in low dimensional model checking that local smoothing tests typically are more sensitive to high-frequency models and empirical process-based tests are more powerful for low-frequency models. Interestingly, our hybrid test $Hybrid_{CF}$ inherits the merits of both local smoothing tests and empirical process-based tests, as it maintains high power for both low-frequency and high-frequency alternatives.

\begin{table}[ht!]\caption{Empirical sizes and powers of the tests $TCvM^{2}_C$, $TCvM^{2}_{CF}$, $\hat{T}_{\rm{Fisher}}^C$, $Hybrid_{CF}$, $RP_n$, and $GRP_n$ for $H_{11}$ in Study 1.}\label{table-H11}
\centering
{\small\scriptsize\hspace{8cm}
\renewcommand{\arraystretch}{0.6}\tabcolsep 0.4cm
\begin{tabular}{*{20}{c}}
\hline
&\multicolumn{1}{c}{a}&\multicolumn{1}{c}{n=300}&\multicolumn{1}{c}{n=300}&\multicolumn{1}{c}{n=300}&\multicolumn{1}{c}{n=300} &\multicolumn{1}{c}{n=300} &\multicolumn{1}{c}{n=300}\\
&&\multicolumn{1}{c}{p=50}&\multicolumn{1}{c}{p=100}&\multicolumn{1}{c}{p=300}&\multicolumn{1}{c}{p=600} &\multicolumn{1}{c}{p=900} &\multicolumn{1}{c}{p=1200}  \\
\hline
$TCvM^{2}_C,\ \Sigma = I_p$
& 0.0 & 0.052 & 0.040 & 0.042 & 0.051 & 0.046 & 0.047 \\
& 0.5 & 0.664 & 0.581 & 0.401 & 0.357 & 0.334 & 0.279 \\
& 1.0 & 0.991 & 0.953 & 0.831 & 0.716 & 0.628 & 0.608 \\
\\
$TCvM^{2}_{CF},\ \Sigma = I_p$
& 0.0 & 0.047 & 0.035 & 0.035 & 0.056 & 0.048 & 0.049 \\
& 0.5 & 0.847 & 0.731 & 0.547 & 0.474 & 0.432 & 0.396 \\
& 1.0 & 1.000 & 0.995 & 0.929 & 0.875 & 0.791 & 0.773 \\
\\
$\hat{T}_{\rm{Fisher}}^C,\ \Sigma = I_p$
& 0.0 & 0.021 & 0.034 & 0.026 & 0.036 & 0.026 & 0.034 \\
& 0.5 & 0.195 & 0.143 & 0.111 & 0.104 & 0.084 & 0.073 \\
& 1.0 & 0.917 & 0.791 & 0.534 & 0.419 & 0.343 & 0.300 \\
\\
$Hybrid_{CF},\ \Sigma = I_p$
& 0.0 & 0.057 & 0.050 & 0.058 & 0.065 & 0.054 & 0.063 \\
& 0.5 & 0.772 & 0.679 & 0.478 & 0.405 & 0.362 & 0.336 \\
& 1.0 & 0.999 & 0.993 & 0.917 & 0.828 & 0.744 & 0.730 \\
\\
$RP_n,\ \Sigma = I_p$
& 0.0 & 0.038 & 0.039 & 0.034 & 0.024 & 0.028 & 0.041 \\
& 0.5 & 0.094 & 0.086 & 0.058 & 0.054 & 0.064 & 0.057 \\
& 1.0 & 0.176 & 0.149 & 0.144 & 0.127 & 0.110 & 0.114 \\
\\
$GRP_n,\ \Sigma = I_p$
& 0.0 & 0.034 & 0.020 & 0.018 & 0.017 & 0.014 & 0.012 \\
& 0.5 & 0.085 & 0.049 & 0.025 & 0.019 & 0.018 & 0.025 \\
& 1.0 & 0.291 & 0.182 & 0.079 & 0.054 & 0.040 & 0.055 \\
\hline
$TCvM^{2}_{C},\ \Sigma = (0.4^{|i-j|})_{p\times p}$
& 0.0 & 0.039 & 0.051 & 0.045 & 0.044 & 0.030 & 0.037 \\
& 0.5 & 0.994 & 0.984 & 0.903 & 0.857 & 0.821 & 0.808 \\
& 1.0 & 1.000 & 1.000 & 0.982 & 0.951 & 0.924 & 0.904 \\
\\
$TCvM^{2}_{CF},\ \Sigma =  (0.4^{|i-j|})_{p\times p}$
& 0.0 & 0.043 & 0.048 & 0.046 & 0.046 & 0.037 & 0.037 \\
& 0.5 & 0.999 & 1.000 & 0.981 & 0.959 & 0.947 & 0.937 \\
& 1.0 & 1.000 & 1.000 & 0.997 & 0.990 & 0.982 & 0.981 \\
\\
$\hat{T}_{\rm{Fisher}}^C,\ \Sigma =   (0.4^{|i-j|})_{p\times p}$
& 0.0 & 0.025 & 0.029 & 0.031 & 0.026 & 0.033 & 0.034 \\
& 0.5 & 0.896 & 0.779 & 0.627 & 0.574 & 0.529 & 0.487 \\
& 1.0 & 1.000 & 1.000 & 0.985 & 0.918 & 0.887 & 0.874 \\
\\
$Hybrid_{CF},\ \Sigma = (0.4^{|i-j|})_{p\times p}$
& 0.0 & 0.055 & 0.058 & 0.057 & 0.053 & 0.047 & 0.050 \\
& 0.5 & 0.999 & 1.000 & 0.975 & 0.950 & 0.929 & 0.908 \\
& 1.0 & 1.000 & 1.000 & 0.996 & 0.983 & 0.980 & 0.978 \\
\\
$RP_n,\ \Sigma =  (0.4^{|i-j|})_{p\times p}$
& 0.0 & 0.029 & 0.037 & 0.035 & 0.037 & 0.039 & 0.042 \\
& 0.5 & 0.470 & 0.391 & 0.340 & 0.314 & 0.327 & 0.318 \\
& 1.0 & 0.601 & 0.534 & 0.481 & 0.436 & 0.445 & 0.485 \\
\\
$GRP_n,\ \Sigma =  (0.4^{|i-j|})_{p\times p}$
& 0.0 & 0.026 & 0.031 & 0.030 & 0.025 & 0.019 & 0.019 \\
& 0.5 & 0.785 & 0.604 & 0.280 & 0.162 & 0.139 & 0.121 \\
& 1.0 & 0.998 & 0.981 & 0.700 & 0.582 & 0.542 & 0.475 \\
\hline
$TCvM^{2}_{C},\ \Sigma = (0.8^{|i-j|})_{p\times p}$
& 0.0 & 0.047 & 0.050 & 0.039 & 0.030 & 0.048 & 0.047 \\
& 0.5 & 1.000 & 1.000 & 0.998 & 0.997 & 0.988 & 0.977 \\
& 1.0 & 1.000 & 1.000 & 1.000 & 0.994 & 0.990 & 0.976 \\
\\
$TCvM^{2}_{CF},\ \Sigma = (0.8^{|i-j|})_{p\times p}$
& 0.0 & 0.040 & 0.046 & 0.037 & 0.035 & 0.056 & 0.049 \\
& 0.5 & 1.000 & 1.000 & 1.000 & 1.000 & 0.999 & 1.000 \\
& 1.0 & 1.000 & 1.000 & 1.000 & 0.999 & 0.999 & 0.997 \\
\\
$\hat{T}_{\rm{Fisher}}^C,\ \Sigma =   (0.8^{|i-j|})_{p\times p}$
& 0.0 & 0.028 & 0.028 & 0.025 & 0.032 & 0.041 & 0.038 \\
& 0.5 & 1.000 & 1.000 & 0.997 & 0.989 & 0.970 & 0.968 \\
& 1.0 & 1.000 & 1.000 & 1.000 & 0.998 & 0.992 & 0.984 \\
\\
$Hybrid_{CF},\ \Sigma = (0.8^{|i-j|})_{p\times p}$
& 0.0 & 0.050 & 0.055 & 0.040 & 0.055 & 0.068 & 0.069 \\
& 0.5 & 1.000 & 1.000 & 1.000 & 1.000 & 0.998 & 0.998 \\
& 1.0 & 1.000 & 1.000 & 1.000 & 0.999 & 0.999 & 0.998 \\
\\
$RP_n,\ \Sigma = (0.8^{|i-j|})_{p\times p}$
& 0.0 & 0.039 & 0.029 & 0.039 & 0.032 & 0.035 & 0.034 \\
& 0.5 & 0.746 & 0.654 & 0.578 & 0.574 & 0.582 & 0.541 \\
& 1.0 & 0.741 & 0.653 & 0.603 & 0.561 & 0.547 & 0.552 \\
\\
$GRP_n,\ \Sigma = (0.8^{|i-j|})_{p\times p}$
& 0.0 & 0.043 & 0.041 & 0.027 & 0.024 & 0.025 & 0.037 \\
& 0.5 & 1.000 & 1.000 & 0.987 & 0.969 & 0.961 & 0.950 \\
& 1.0 & 1.000 & 1.000 & 1.000 & 0.997 & 0.995 & 0.999 \\
\hline
\end{tabular}}
\end{table}

\begin{table}[ht!]\caption{Empirical sizes and powers of the tests $TCvM^{2}_C$, $TCvM^{2}_{CF}$, $\hat{T}_{\rm{Fisher}}^C$, $Hybrid_{CF}$, $RP_n$ and $GRP_n$ for $H_{12}$ in Study 1.}\label{table-H12}
\centering
{\small\scriptsize\hspace{8cm}
\renewcommand{\arraystretch}{0.6}\tabcolsep 0.4cm
\begin{tabular}{*{20}{c}}
\hline
&\multicolumn{1}{c}{a}&\multicolumn{1}{c}{n=300}&\multicolumn{1}{c}{n=300}&\multicolumn{1}{c}{n=300}&\multicolumn{1}{c}{n=300} &\multicolumn{1}{c}{n=300} &\multicolumn{1}{c}{n=300}\\
&&\multicolumn{1}{c}{p=50}&\multicolumn{1}{c}{p=100}&\multicolumn{1}{c}{p=300}&\multicolumn{1}{c}{p=600} &\multicolumn{1}{c}{p=900} &\multicolumn{1}{c}{p=1200}  \\
\hline
$TCvM^{2}_{C},\ \Sigma = I_p$
& 0.0 & 0.033 & 0.045 & 0.052 & 0.042 & 0.048 & 0.031 \\
& 0.5 & 0.055 & 0.052 & 0.039 & 0.033 & 0.046 & 0.061 \\
& 1.0 & 0.046 & 0.059 & 0.035 & 0.044 & 0.043 & 0.035 \\
\\
$TCvM^{2}_{CF},\ \Sigma = I_p$
& 0.0 & 0.048 & 0.049 & 0.042 & 0.045 & 0.046 & 0.038 \\
& 0.5 & 0.054 & 0.055 & 0.043 & 0.048 & 0.032 & 0.051 \\
& 1.0 & 0.052 & 0.052 & 0.037 & 0.051 & 0.048 & 0.043 \\
\\
$\hat{T}_{\rm{Fisher}}^C,\ \Sigma = I_p$
& 0.0 & 0.020 & 0.017 & 0.031 & 0.036 & 0.037 & 0.037 \\
& 0.5 & 0.320 & 0.216 & 0.085 & 0.074 & 0.053 & 0.053 \\
& 1.0 & 0.838 & 0.650 & 0.227 & 0.124 & 0.091 & 0.095 \\
\\
$Hybrid_{CF},\ \Sigma = I_p$
& 0.0 & 0.050 & 0.055 & 0.050 & 0.059 & 0.059 & 0.050 \\
& 0.5 & 0.274 & 0.182 & 0.096 & 0.075 & 0.052 & 0.059 \\
& 1.0 & 0.689 & 0.413 & 0.178 & 0.116 & 0.083 & 0.085 \\
\\
$RP_n,\ \Sigma = I_p$
& 0.0 & 0.036 & 0.030 & 0.039 & 0.040 & 0.038 & 0.032 \\
& 0.5 & 0.045 & 0.037 & 0.032 & 0.042 & 0.029 & 0.040 \\
& 1.0 & 0.033 & 0.058 & 0.032 & 0.049 & 0.042 & 0.036 \\
\\
$GRP_n,\ \Sigma = I_p$
& 0.0 & 0.037 & 0.033 & 0.018 & 0.010 & 0.016 & 0.012 \\
& 0.5 & 0.040 & 0.034 & 0.027 & 0.019 & 0.014 & 0.015 \\
& 1.0 & 0.069 & 0.066 & 0.045 & 0.027 & 0.032 & 0.027 \\
\hline
$TCvM^{2}_{C},\ \Sigma = (0.4^{|i-j|})_{p\times p}$
& 0.0 & 0.048 & 0.052 & 0.041 & 0.037 & 0.041 & 0.058 \\
& 0.5 & 0.042 & 0.043 & 0.039 & 0.038 & 0.049 & 0.057 \\
& 1.0 & 0.052 & 0.046 & 0.037 & 0.037 & 0.039 & 0.043 \\
\\
$TCvM^{2}_{CF},\ \Sigma =  (0.4^{|i-j|})_{p\times p}$
& 0.0 & 0.044 & 0.053 & 0.047 & 0.038 & 0.044 & 0.040 \\
& 0.5 & 0.037 & 0.043 & 0.043 & 0.042 & 0.037 & 0.057 \\
& 1.0 & 0.052 & 0.055 & 0.044 & 0.042 & 0.035 & 0.043 \\
\\
$\hat{T}_{\rm{Fisher}}^C,\ \Sigma =   (0.4^{|i-j|})_{p\times p}$
& 0.0 & 0.031 & 0.027 & 0.030 & 0.031 & 0.037 & 0.037 \\
& 0.5 & 0.350 & 0.279 & 0.155 & 0.142 & 0.114 & 0.106 \\
& 1.0 & 0.924 & 0.837 & 0.567 & 0.404 & 0.308 & 0.258 \\
\\
$Hybrid_{CF},\ \Sigma = (0.4^{|i-j|})_{p\times p}$
& 0.0 & 0.054 & 0.060 & 0.058 & 0.059 & 0.058 & 0.068 \\
& 0.5 & 0.289 & 0.197 & 0.123 & 0.121 & 0.094 & 0.094 \\
& 1.0 & 0.753 & 0.563 & 0.315 & 0.210 & 0.160 & 0.153 \\
\\
$RP_n,\ \Sigma =  (0.4^{|i-j|})_{p\times p}$
& 0.0 & 0.036 & 0.039 & 0.033 & 0.041 & 0.042 & 0.041 \\
& 0.5 & 0.040 & 0.051 & 0.019 & 0.039 & 0.041 & 0.045 \\
& 1.0 & 0.035 & 0.037 & 0.037 & 0.041 & 0.048 & 0.041 \\
\\
$GRP_n,\ \Sigma =  (0.4^{|i-j|})_{p\times p}$
& 0.0 & 0.033 & 0.034 & 0.021 & 0.030 & 0.022 & 0.016 \\
& 0.5 & 0.043 & 0.053 & 0.032 & 0.018 & 0.025 & 0.021 \\
& 1.0 & 0.075 & 0.056 & 0.062 & 0.036 & 0.027 & 0.038 \\
\hline
$TCvM^{2}_{C},\ \Sigma = (0.8^{|i-j|})_{p\times p}$
& 0.0 & 0.032 & 0.047 & 0.049 & 0.038 & 0.043 & 0.045 \\
& 0.5 & 0.042 & 0.053 & 0.042 & 0.048 & 0.059 & 0.052 \\
& 1.0 & 0.051 & 0.048 & 0.039 & 0.040 & 0.045 & 0.034 \\
\\
$TCvM^{2}_{CF},\ \Sigma = (0.8^{|i-j|})_{p\times p}$
& 0.0 & 0.036 & 0.046 & 0.040 & 0.041 & 0.054 & 0.045 \\
& 0.5 & 0.037 & 0.046 & 0.046 & 0.049 & 0.042 & 0.041 \\
& 1.0 & 0.047 & 0.044 & 0.031 & 0.052 & 0.049 & 0.032 \\
\\
$\hat{T}_{\rm{Fisher}}^C,\ \Sigma =   (0.8^{|i-j|})_{p\times p}$
& 0.0 & 0.024 & 0.033 & 0.032 & 0.034 & 0.028 & 0.040 \\
& 0.5 & 0.307 & 0.247 & 0.246 & 0.178 & 0.150 & 0.169 \\
& 1.0 & 0.911 & 0.872 & 0.773 & 0.693 & 0.633 & 0.581 \\
\\
$Hybrid_{CF},\ \Sigma = (0.8^{|i-j|})_{p\times p}$
& 0.0 & 0.054 & 0.056 & 0.066 & 0.057 & 0.051 & 0.069 \\
& 0.5 & 0.332 & 0.267 & 0.216 & 0.153 & 0.135 & 0.123 \\
& 1.0 & 0.881 & 0.767 & 0.536 & 0.399 & 0.403 & 0.306 \\
\\
$RP_n,\ \Sigma = (0.8^{|i-j|})_{p\times p}$
& 0.0 & 0.035 & 0.052 & 0.051 & 0.049 & 0.046 & 0.035 \\
& 0.5 & 0.050 & 0.029 & 0.038 & 0.042 & 0.037 & 0.032 \\
& 1.0 & 0.055 & 0.034 & 0.036 & 0.044 & 0.040 & 0.033 \\
\\
$GRP_n,\ \Sigma = (0.8^{|i-j|})_{p\times p}$
& 0.0 & 0.038 & 0.032 & 0.032 & 0.034 & 0.027 & 0.035 \\
& 0.5 & 0.065 & 0.039 & 0.047 & 0.038 & 0.037 & 0.043 \\
& 1.0 & 0.074 & 0.068 & 0.061 & 0.064 & 0.059 & 0.057 \\
\hline
\end{tabular}}
\end{table}

\begin{table}[ht!]\caption{Empirical sizes and powers of the tests $TCvM^{2}_C$, $TCvM^{2}_{CF}$, $\hat{T}_{\rm{Fisher}}^C$, $Hybrid_{CF}$, $RP_n$ and $GRP_n$ for $H_{13}$ in Study 1.}\label{table-H13}
\centering
{\small\scriptsize\hspace{8cm}
\renewcommand{\arraystretch}{0.6}\tabcolsep 0.4cm
\begin{tabular}{*{20}{c}}
\hline
&\multicolumn{1}{c}{a}&\multicolumn{1}{c}{n=300}&\multicolumn{1}{c}{n=300}&\multicolumn{1}{c}{n=300}&\multicolumn{1}{c}{n=300} &\multicolumn{1}{c}{n=300} &\multicolumn{1}{c}{n=300}\\
&&\multicolumn{1}{c}{p=50}&\multicolumn{1}{c}{p=100}&\multicolumn{1}{c}{p=300}&\multicolumn{1}{c}{p=600} &\multicolumn{1}{c}{p=900} &\multicolumn{1}{c}{p=1200}  \\
\hline
$TCvM^{2}_{C},\ \Sigma = I_p$
& 0.0 & 0.047 & 0.043 & 0.046 & 0.043 & 0.038 & 0.040 \\
& 0.5 & 0.312 & 0.211 & 0.135 & 0.150 & 0.115 & 0.104 \\
& 1.0 & 0.776 & 0.668 & 0.401 & 0.310 & 0.267 & 0.243 \\
\\
$TCvM^{2}_{CF},\ \Sigma = I_p$
& 0.0 & 0.061 & 0.037 & 0.045 & 0.039 & 0.037 & 0.045 \\
& 0.5 & 0.397 & 0.294 & 0.189 & 0.177 & 0.132 & 0.120 \\
& 1.0 & 0.920 & 0.834 & 0.545 & 0.417 & 0.387 & 0.310 \\
\\
$\hat{T}_{\rm{Fisher}}^C,\ \Sigma = I_p$
& 0.0 & 0.025 & 0.021 & 0.025 & 0.022 & 0.025 & 0.031 \\
& 0.5 & 0.056 & 0.038 & 0.037 & 0.038 & 0.026 & 0.029 \\
& 1.0 & 0.189 & 0.095 & 0.061 & 0.063 & 0.051 & 0.050 \\
\\
$Hybrid_{CF},\ \Sigma = I_p$
& 0.0 & 0.055 & 0.051 & 0.061 & 0.044 & 0.044 & 0.052 \\
& 0.5 & 0.334 & 0.226 & 0.164 & 0.136 & 0.111 & 0.110 \\
& 1.0 & 0.892 & 0.778 & 0.483 & 0.365 & 0.332 & 0.258 \\
\\
$RP_n,\ \Sigma = I_p$
& 0.0 & 0.043 & 0.044 & 0.033 & 0.044 & 0.034 & 0.048 \\
& 0.5 & 0.054 & 0.045 & 0.061 & 0.046 & 0.042 & 0.050 \\
& 1.0 & 0.072 & 0.061 & 0.069 & 0.058 & 0.065 & 0.053 \\
\\
$GRP_n,\ \Sigma = I_p$
& 0.0 & 0.026 & 0.048 & 0.019 & 0.010 & 0.008 & 0.014 \\
& 0.5 & 0.065 & 0.056 & 0.022 & 0.024 & 0.016 & 0.018 \\
& 1.0 & 0.108 & 0.097 & 0.040 & 0.040 & 0.028 & 0.029 \\
\hline
$TCvM^{2}_{C},\ \Sigma = (0.4^{|i-j|})_{p\times p}$
& 0.0 & 0.037 & 0.043 & 0.052 & 0.045 & 0.043 & 0.042 \\
& 0.5 & 0.943 & 0.866 & 0.741 & 0.647 & 0.613 & 0.574 \\
& 1.0 & 0.998 & 0.996 & 0.943 & 0.869 & 0.798 & 0.803 \\
\\
$TCvM^{2}_{CF},\ \Sigma =  (0.4^{|i-j|})_{p\times p}$
& 0.0 & 0.038 & 0.040 & 0.045 & 0.038 & 0.041 & 0.047 \\
& 0.5 & 0.990 & 0.976 & 0.895 & 0.820 & 0.790 & 0.758 \\
& 1.0 & 1.000 & 1.000 & 0.992 & 0.974 & 0.946 & 0.945 \\
\\
$\hat{T}_{\rm{Fisher}}^C,\ \Sigma =   (0.4^{|i-j|})_{p\times p}$
& 0.0 & 0.031 & 0.026 & 0.028 & 0.033 & 0.030 & 0.037 \\
& 0.5 & 0.437 & 0.306 & 0.187 & 0.143 & 0.129 & 0.138 \\
& 1.0 & 0.954 & 0.831 & 0.586 & 0.409 & 0.350 & 0.305 \\
\\
$Hybrid_{CF},\ \Sigma = (0.4^{|i-j|})_{p\times p}$
& 0.0 & 0.044 & 0.060 & 0.061 & 0.057 & 0.056 & 0.056 \\
& 0.5 & 0.981 & 0.959 & 0.855 & 0.781 & 0.735 & 0.716 \\
& 1.0 & 1.000 & 1.000 & 0.986 & 0.960 & 0.924 & 0.918 \\
\\
$RP_n,\ \Sigma =  (0.4^{|i-j|})_{p\times p}$
& 0.0 & 0.037 & 0.039 & 0.037 & 0.042 & 0.041 & 0.033 \\
& 0.5 & 0.226 & 0.214 & 0.160 & 0.166 & 0.159 & 0.179 \\
& 1.0 & 0.284 & 0.237 & 0.214 & 0.203 & 0.236 & 0.228 \\
\\
$GRP_n,\ \Sigma =  (0.4^{|i-j|})_{p\times p}$
& 0.0 & 0.030 & 0.030 & 0.015 & 0.021 & 0.020 & 0.021 \\
& 0.5 & 0.642 & 0.462 & 0.211 & 0.155 & 0.134 & 0.129 \\
& 1.0 & 0.933 & 0.800 & 0.430 & 0.312 & 0.320 & 0.269 \\
\hline
$TCvM^{2}_{C},\ \Sigma = (0.8^{|i-j|})_{p\times p}$
& 0.0 & 0.041 & 0.035 & 0.042 & 0.039 & 0.040 & 0.037 \\
& 0.5 & 1.000 & 0.999 & 0.989 & 0.964 & 0.954 & 0.941 \\
& 1.0 & 1.000 & 0.999 & 0.987 & 0.954 & 0.934 & 0.921 \\
\\
$TCvM^{2}_{CF},\ \Sigma = (0.8^{|i-j|})_{p\times p}$
& 0.0 & 0.048 & 0.048 & 0.044 & 0.044 & 0.040 & 0.034 \\
& 0.5 & 1.000 & 1.000 & 1.000 & 0.999 & 0.991 & 0.996 \\
& 1.0 & 1.000 & 1.000 & 0.998 & 0.997 & 0.990 & 0.986 \\
\\
$\hat{T}_{\rm{Fisher}}^C,\ \Sigma =   (0.8^{|i-j|})_{p\times p}$
& 0.0 & 0.036 & 0.024 & 0.034 & 0.032 & 0.034 & 0.042 \\
& 0.5 & 1.000 & 1.000 & 0.991 & 0.969 & 0.928 & 0.887 \\
& 1.0 & 1.000 & 1.000 & 0.998 & 0.983 & 0.962 & 0.935 \\
\\
$Hybrid_{CF},\ \Sigma = (0.8^{|i-j|})_{p\times p}$
& 0.0 & 0.053 & 0.052 & 0.047 & 0.053 & 0.066 & 0.057 \\
& 0.5 & 1.000 & 1.000 & 1.000 & 0.998 & 0.990 & 0.990 \\
& 1.0 & 1.000 & 1.000 & 0.998 & 0.998 & 0.988 & 0.978 \\
\\
$RP_n,\ \Sigma = (0.8^{|i-j|})_{p\times p}$
& 0.0 & 0.043 & 0.042 & 0.048 & 0.032 & 0.038 & 0.036 \\
& 0.5 & 0.663 & 0.622 & 0.513 & 0.506 & 0.521 & 0.482 \\
& 1.0 & 0.687 & 0.603 & 0.557 & 0.523 & 0.525 & 0.531 \\
\\
$GRP_n,\ \Sigma = (0.8^{|i-j|})_{p\times p}$
& 0.0 & 0.053 & 0.037 & 0.032 & 0.024 & 0.029 & 0.032 \\
& 0.5 & 1.000 & 1.000 & 0.993 & 0.985 & 0.974 & 0.975 \\
& 1.0 & 1.000 & 1.000 & 1.000 & 0.998 & 0.996 & 0.989 \\
\hline
\end{tabular}}
\end{table}

Next, we investigate the finite sample performance of our proposed tests for the goodness-of-fit of logistic regression models.

{\em Study 2.} The data are generated from the logistic regression model according to
$$ Y|X \sim Bernoulli\{ \mu(\beta_0^{\top}X + ag(X)) \},$$
where $\mu(z) = 1/(1+\exp(-z))$. We consider two different cases for the misspecified $g(X)$:
\begin{eqnarray*}
H_{21}:  g(X) &=& 0.2(\beta_0^{\top}X)^2,   \\
%H_{22}:  g(X) &=& 2\cos(0.6 \pi \beta_0^{\top}X),
H_{22}:  g(X) &=& X^{(1)} X^{(2)} + X^{(2)}X^{(3)} + X^{(3)}X^{(4)} + X^{(4)}X^{(5)},
\end{eqnarray*}
where the parameter $\beta_0 = (1,1,1,1,1,0,\dots,0)^{\top}$ and the covariates $X$ are the same as in study 1, and the sample size $n=600$ with the covariate dimension $p \in \{50, 100, 300, 600, 900, 1200\}$.

Since the $RP_n$ test proposed by \cite{Shah2018} cannot be applied in testing GLMs, we only compare our tests with the tests $GRP_n$ and $\hat{T}_{\rm{Fisher}}^C$ developed by \cite{tan2025b} and \cite{Jankova2020}, respectively. The simulation results are provided in Tables 4-5. It can be seen that our tests $TCvM^{2}_C$, $TCvM^{2}_{CF}$, $Hybrid_{CF}$, and the local smoothing test $\hat{T}_{\rm{Fisher}}^C$ control the empirical size in most cases. The tests $TCvM^{2}_C$ and $TCvM^{2}_{CF}$ are generally conservative, exhibiting smaller empirical sizes in large dimensional settings. The test $\hat{T}_{\rm{Fisher}}^C$ becomes liberal with large empirical sizes in settings with high correlation ($\rho = 0.8$) and large dimension ($p=1200$). However, the empirical size of $GRP_n$ is far from the significant level when the covariate correlation $\rho = 0.4 $ or $\rho = 0.8$. For the empirical power, our tests $TCvM^{2}_{CF}$ and $Hybrid_{CF}$ typically have higher power than the other competitors. Moreover, the empirical powers of all tests increase as the correlation of the covariates $X$ grows.

\begin{table}[ht!]\caption{Empirical sizes and powers of the tests $TCvM^{2}_C$, $TCvM^{2}_{CF}$, $\hat{T}_{\rm{Fisher}}^C$, $Hybrid_{CF}$, and $GRP_n$ for $H_{21}$ in Study 2.}\label{table-H21}
\centering
{\small\scriptsize\hspace{8cm}
\renewcommand{\arraystretch}{0.6}\tabcolsep 0.4cm
\begin{tabular}{*{20}{c}}
\hline
&\multicolumn{1}{c}{a}&\multicolumn{1}{c}{n=600}&\multicolumn{1}{c}{n=600}&\multicolumn{1}{c}{n=600}&\multicolumn{1}{c}{n=600} &\multicolumn{1}{c}{n=600} &\multicolumn{1}{c}{n=600}\\
&&\multicolumn{1}{c}{p=50}&\multicolumn{1}{c}{p=100}&\multicolumn{1}{c}{p=300}&\multicolumn{1}{c}{p=600} &\multicolumn{1}{c}{p=900} &\multicolumn{1}{c}{p=1200}  \\
\hline
$TCvM^{2}_{C},\ \Sigma = I_p$
& 0.0 & 0.068 & 0.065 & 0.033 & 0.013 & 0.029 & 0.027 \\
& 0.5 & 0.296 & 0.209 & 0.084 & 0.043 & 0.038 & 0.036 \\
& 1.0 & 0.800 & 0.630 & 0.338 & 0.226 & 0.156 & 0.155 \\
\\
$TCvM^{2}_{CF},\ \Sigma = I_p$
& 0.0 & 0.087 & 0.065 & 0.016 & 0.011 & 0.009 & 0.011 \\
& 0.5 & 0.348 & 0.264 & 0.085 & 0.045 & 0.024 & 0.019 \\
& 1.0 & 0.944 & 0.816 & 0.493 & 0.266 & 0.165 & 0.153 \\
\\
$\hat{T}_{\rm{Fisher}}^C,\ \Sigma = I_p$
& 0.0 & 0.028 & 0.028 & 0.026 & 0.032 & 0.040 & 0.049 \\
& 0.5 & 0.084 & 0.077 & 0.032 & 0.034 & 0.034 & 0.045 \\
& 1.0 & 0.681 & 0.489 & 0.203 & 0.101 & 0.068 & 0.058 \\
\\
$Hybrid_{CF},\ \Sigma = I_p$
& 0.0 & 0.080 & 0.053 & 0.041 & 0.024 & 0.028 & 0.028 \\
& 0.5 & 0.278 & 0.236 & 0.088 & 0.066 & 0.041 & 0.033 \\
& 1.0 & 0.924 & 0.787 & 0.445 & 0.257 & 0.169 & 0.144 \\
\\
$GRP_n,\ \Sigma = I_p$
& 0.0 & 0.084 & 0.065 & 0.030 & 0.030 & 0.018 & 0.020 \\
& 0.5 & 0.068 & 0.048 & 0.034 & 0.027 & 0.027 & 0.023 \\
& 1.0 & 0.077 & 0.042 & 0.029 & 0.025 & 0.024 & 0.024 \\
\hline
$TCvM^{2}_{C},\ \Sigma = (0.4^{|i-j|})_{p\times p}$
& 0.0 & 0.093 & 0.076 & 0.042 & 0.023 & 0.011 & 0.011 \\
& 0.5 & 0.566 & 0.444 & 0.259 & 0.161 & 0.102 & 0.085 \\
& 1.0 & 0.960 & 0.902 & 0.711 & 0.537 & 0.449 & 0.408 \\
\\
$TCvM^{2}_{CF},\ \Sigma =  (0.4^{|i-j|})_{p\times p}$
& 0.0 & 0.088 & 0.050 & 0.033 & 0.006 & 0.008 & 0.005 \\
& 0.5 & 0.673 & 0.562 & 0.313 & 0.174 & 0.102 & 0.075 \\
& 1.0 & 1.000 & 0.992 & 0.950 & 0.776 & 0.677 & 0.581 \\
\\
$\hat{T}_{\rm{Fisher}}^C,\ \Sigma =   (0.4^{|i-j|})_{p\times p}$
& 0.0 & 0.027 & 0.023 & 0.033 & 0.030 & 0.056 & 0.076 \\
& 0.5 & 0.215 & 0.224 & 0.146 & 0.090 & 0.076 & 0.090 \\
& 1.0 & 0.996 & 0.961 & 0.810 & 0.682 & 0.588 & 0.519 \\
\\
$Hybrid_{CF},\ \Sigma = (0.4^{|i-j|})_{p\times p}$
& 0.0 & 0.087 & 0.054 & 0.049 & 0.038 & 0.037 & 0.034 \\
& 0.5 & 0.620 & 0.503 & 0.298 & 0.187 & 0.112 & 0.101 \\
& 1.0 & 0.999 & 0.991 & 0.948 & 0.809 & 0.714 & 0.630 \\
\\
$GRP_n,\ \Sigma =  (0.4^{|i-j|})_{p\times p}$
& 0.0 & 0.225 & 0.168 & 0.114 & 0.105 & 0.095 & 0.078 \\
& 0.5 & 0.151 & 0.103 & 0.057 & 0.036 & 0.036 & 0.038 \\
& 1.0 & 0.776 & 0.478 & 0.167 & 0.104 & 0.077 & 0.072 \\
\hline
$TCvM^{2}_{C},\ \Sigma = (0.8^{|i-j|})_{p\times p}$
& 0.0 & 0.088 & 0.073 & 0.051 & 0.037 & 0.031 & 0.028 \\
& 0.5 & 0.924 & 0.923 & 0.833 & 0.757 & 0.644 & 0.591 \\
& 1.0 & 0.980 & 0.965 & 0.860 & 0.771 & 0.724 & 0.664 \\
\\
$TCvM^{2}_{CF},\ \Sigma = (0.8^{|i-j|})_{p\times p}$
& 0.0 & 0.088 & 0.055 & 0.031 & 0.018 & 0.012 & 0.008 \\
& 0.5 & 0.991 & 0.981 & 0.950 & 0.877 & 0.776 & 0.702 \\
& 1.0 & 1.000 & 1.000 & 0.988 & 0.971 & 0.946 & 0.914 \\
\\
$\hat{T}_{\rm{Fisher}}^C,\ \Sigma =   (0.8^{|i-j|})_{p\times p}$
& 0.0 & 0.034 & 0.031 & 0.026 & 0.058 & 0.079 & 0.130 \\
& 0.5 & 0.935 & 0.921 & 0.874 & 0.872 & 0.823 & 0.757 \\
& 1.0 & 1.000 & 0.998 & 0.954 & 0.910 & 0.879 & 0.854 \\
\\
$Hybrid_{CF},\ \Sigma = (0.8^{|i-j|})_{p\times p}$
& 0.0 & 0.081 & 0.058 & 0.035 & 0.048 & 0.049 & 0.063 \\
& 0.5 & 0.991 & 0.983 & 0.958 & 0.912 & 0.858 & 0.801 \\
& 1.0 & 1.000 & 0.999 & 0.988 & 0.971 & 0.957 & 0.933 \\
\\
$GRP_n,\ \Sigma = (0.8^{|i-j|})_{p\times p}$
& 0.0 & 0.410 & 0.393 & 0.357 & 0.322 & 0.310 & 0.305 \\
& 0.5 & 0.830 & 0.769 & 0.606 & 0.461 & 0.366 & 0.305 \\
& 1.0 & 1.000 & 1.000 & 0.964 & 0.970 & 0.938 & 0.893 \\
\hline
\end{tabular}}
\end{table}

\begin{table}[ht!]\caption{Empirical sizes and powers of the tests $TCvM^{2}_C$, $TCvM^{2}_{CF}$, $\hat{T}_{\rm{Fisher}}^C$, $Hybrid_{CF}$, and $GRP_n$ for $H_{22}$ in Study 2.}\label{table-H22}
\centering
{\small\scriptsize\hspace{8cm}
\renewcommand{\arraystretch}{0.6}\tabcolsep 0.4cm
\begin{tabular}{*{20}{c}}
\hline
&\multicolumn{1}{c}{a}&\multicolumn{1}{c}{n=600}&\multicolumn{1}{c}{n=600}&\multicolumn{1}{c}{n=600}&\multicolumn{1}{c}{n=600} &\multicolumn{1}{c}{n=600} &\multicolumn{1}{c}{n=600}\\
&&\multicolumn{1}{c}{p=50}&\multicolumn{1}{c}{p=100}&\multicolumn{1}{c}{p=300}&\multicolumn{1}{c}{p=600} &\multicolumn{1}{c}{p=900} &\multicolumn{1}{c}{p=1200}  \\
\hline
$TCvM^{2}_{C},\ \Sigma = I_p$
& 0.0 & 0.090 & 0.073 & 0.026 & 0.023 & 0.019 & 0.016 \\
& 0.5 & 0.169 & 0.129 & 0.064 & 0.045 & 0.031 & 0.025 \\
& 1.0 & 0.433 & 0.357 & 0.141 & 0.087 & 0.079 & 0.070 \\
\\
$TCvM^{2}_{CF},\ \Sigma = I_p$
& 0.0 & 0.093 & 0.066 & 0.011 & 0.017 & 0.005 & 0.007 \\
& 0.5 & 0.182 & 0.121 & 0.048 & 0.031 & 0.021 & 0.011 \\
& 1.0 & 0.536 & 0.425 & 0.155 & 0.094 & 0.074 & 0.063 \\
\\
$\hat{T}_{\rm{Fisher}}^C,\ \Sigma = I_p$
& 0.0 & 0.022 & 0.027 & 0.033 & 0.046 & 0.033 & 0.043 \\
& 0.5 & 0.039 & 0.046 & 0.031 & 0.043 & 0.045 & 0.042 \\
& 1.0 & 0.241 & 0.195 & 0.064 & 0.046 & 0.043 & 0.041 \\
\\
$Hybrid_{CF},\ \Sigma = I_p$
& 0.0 & 0.082 & 0.054 & 0.033 & 0.043 & 0.026 & 0.028 \\
& 0.5 & 0.162 & 0.121 & 0.057 & 0.037 & 0.049 & 0.023 \\
& 1.0 & 0.489 & 0.367 & 0.153 & 0.103 & 0.095 & 0.066 \\
\\
$GRP_n,\ \Sigma = I_p$
& 0.0 & 0.076 & 0.065 & 0.035 & 0.029 & 0.025 & 0.020 \\
& 0.5 & 0.130 & 0.078 & 0.033 & 0.036 & 0.020 & 0.022 \\
& 1.0 & 0.235 & 0.094 & 0.037 & 0.031 & 0.027 & 0.019 \\
\hline
$TCvM^{2}_{C},\ \Sigma = (0.4^{|i-j|})_{p\times p}$
& 0.0 & 0.092 & 0.076 & 0.038 & 0.019 & 0.019 & 0.015 \\
& 0.5 & 0.451 & 0.351 & 0.185 & 0.100 & 0.081 & 0.082 \\
& 1.0 & 0.902 & 0.810 & 0.570 & 0.439 & 0.334 & 0.334 \\
\\
$TCvM^{2D}_{CF},\ \Sigma =  (0.4^{|i-j|})_{p\times p}$
& 0.0 & 0.088 & 0.046 & 0.020 & 0.009 & 0.010 & 0.003 \\
& 0.5 & 0.506 & 0.416 & 0.193 & 0.102 & 0.072 & 0.055 \\
& 1.0 & 0.989 & 0.965 & 0.858 & 0.705 & 0.574 & 0.488 \\
\\
$\hat{T}_{\rm{Fisher}}^C,\ \Sigma =   (0.4^{|i-j|})_{p\times p}$
& 0.0 & 0.023 & 0.027 & 0.034 & 0.049 & 0.061 & 0.090 \\
& 0.5 & 0.154 & 0.136 & 0.090 & 0.067 & 0.079 & 0.063 \\
& 1.0 & 0.906 & 0.835 & 0.627 & 0.487 & 0.384 & 0.327 \\
\\
$Hybrid_{CF},\ \Sigma = (0.4^{|i-j|})_{p\times p}$
& 0.0 & 0.071 & 0.051 & 0.036 & 0.041 & 0.038 & 0.028 \\
& 0.5 & 0.491 & 0.376 & 0.183 & 0.114 & 0.098 & 0.069 \\
& 1.0 & 0.985 & 0.953 & 0.866 & 0.710 & 0.585 & 0.531 \\
\\
$GRP_n,\ \Sigma =  (0.4^{|i-j|})_{p\times p}$
& 0.0 & 0.195 & 0.173 & 0.110 & 0.072 & 0.098 & 0.078 \\
& 0.5 & 0.280 & 0.181 & 0.107 & 0.075 & 0.065 & 0.048 \\
& 1.0 & 0.944 & 0.790 & 0.336 & 0.191 & 0.152 & 0.109 \\
\hline
$TCvM^{2}_{C},\ \Sigma = (0.8^{|i-j|})_{p\times p}$
& 0.0 & 0.115 & 0.098 & 0.054 & 0.056 & 0.029 & 0.022 \\
& 0.5 & 0.778 & 0.699 & 0.582 & 0.475 & 0.386 & 0.344 \\
& 1.0 & 0.985 & 0.969 & 0.908 & 0.791 & 0.720 & 0.694 \\
\\
$TCvM^{2}_{CF},\ \Sigma = (0.8^{|i-j|})_{p\times p}$
& 0.0 & 0.089 & 0.077 & 0.028 & 0.024 & 0.012 & 0.007 \\
& 0.5 & 0.906 & 0.833 & 0.703 & 0.581 & 0.455 & 0.398 \\
& 1.0 & 1.000 & 1.000 & 0.996 & 0.982 & 0.935 & 0.912 \\
\\
$\hat{T}_{\rm{Fisher}}^C,\ \Sigma =   (0.8^{|i-j|})_{p\times p}$
& 0.0 & 0.020 & 0.029 & 0.035 & 0.047 & 0.093 & 0.119 \\
& 0.5 & 0.525 & 0.519 & 0.479 & 0.505 & 0.423 & 0.391 \\
& 1.0 & 1.000 & 0.999 & 0.975 & 0.937 & 0.906 & 0.886 \\
\\
$Hybrid_{CF},\ \Sigma = (0.8^{|i-j|})_{p\times p}$
& 0.0 & 0.082 & 0.068 & 0.034 & 0.037 & 0.061 & 0.063 \\
& 0.5 & 0.891 & 0.824 & 0.699 & 0.632 & 0.521 & 0.479 \\
& 1.0 & 1.000 & 1.000 & 0.994 & 0.983 & 0.958 & 0.946 \\
\\
$GRP_n,\ \Sigma = (0.8^{|i-j|})_{p\times p}$
& 0.0 & 0.409 & 0.392 & 0.368 & 0.357 & 0.345 & 0.323 \\
& 0.5 & 0.566 & 0.499 & 0.341 & 0.282 & 0.228 & 0.195 \\
& 1.0 & 1.000 & 1.000 & 0.962 & 0.996 & 0.981 & 0.944 \\
\hline
\end{tabular}}
\end{table}

\subsection{Real data examples}
\begin{figure}
\centering
\includegraphics[width=0.55\linewidth]{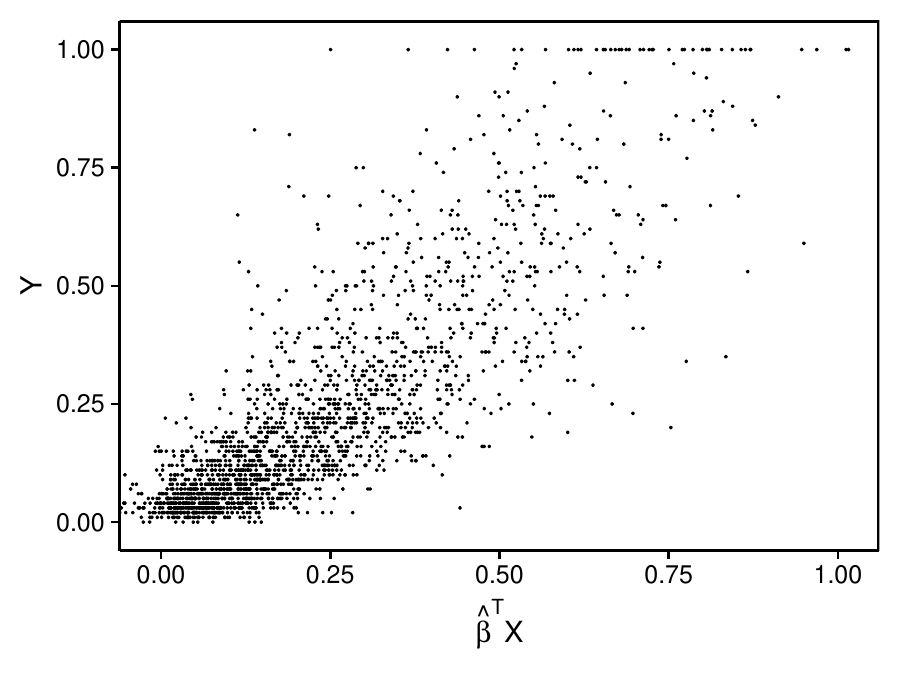}
\caption{The scatter plot of $Y$ versus $\hat{\beta}^{\top}X$.}
\label{fig:y-xbeta}
\vspace{-6pt}
\end{figure}

In this subsection, we evaluate the proposed tests using two real datasets: the Communities and Crime data and the Acute Myeloid Leukemia (AML) data \citep{bottomly2022integrative}.
The Communities and Crime dataset, available from the UCI Machine Learning Repository (\url{https://archive.ics.uci.edu/dataset/183/communities+and+crime}), contains $1994$ observations with one response variable, the per capita violent crime rate, and $99$ predictors describing demographic and law-enforcement characteristics. Let $Y$ denote the per capita violent crime rate and $X=(X^{(1)}, \dots, X^{(99)})^{\top}$ represent the predictor vector. We first evaluate the adequacy of a sparse linear regression model, $Y = a_0 + \beta_0^{\top}X + \varepsilon$, for this dataset using our tests. The choices for the kernel function and bandwidth are the same as those in the simulation studies.
The $p$-values of $TCvM^{2}_{C}$, $TCvM^{2}_{CF}$, and $Hybrid_{CF}$ are approximately $0.0037$, $0.0053$, and $0.0001$, respectively. These results strongly reject the null hypothesis, indicating that the linear relationship between $Y$ and $X$ is not adequate for fitting this dataset. Figure~\ref{fig:y-xbeta} presents a scatter plot of the response variable $Y$ versus $\hat{\beta}_0^{\top}X$, where $\hat{\beta}_0$ is a post-Lasso estimator from the linear model. This plot also suggests that a linear relationship between $Y$ and $X$ may not be plausible. Furthermore, it suggests the potential existence of a quadratic relationship between $Y$ and $\hat{\beta}_0^{\top}X$.

\begin{figure}
\centering
\includegraphics[width=0.55\linewidth]{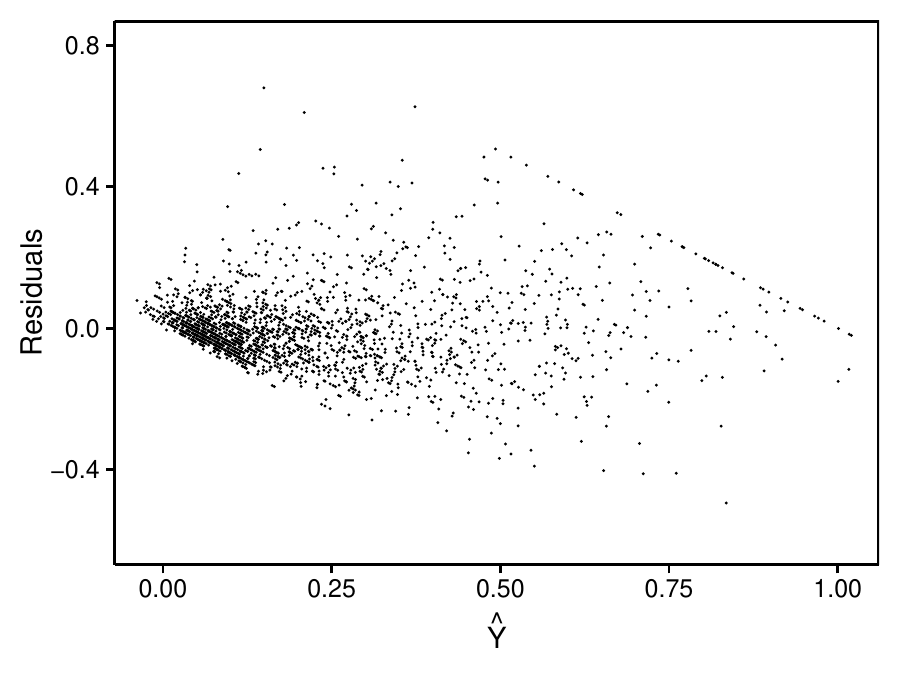}
\caption{The scatter plot of the residuals from the quadratic polynomial model versus the fitted values $\hat{Y}$.}
\label{fig:Res}
\vspace{-6pt}
\end{figure}

To identify a more appropriate relationship between $Y$ and $X$, we expand the model by incorporating the quadratic and interaction terms of covariates, leading to the following quadratic polynomial model
\begin{equation}\label{quadratic-model}
Y = a_0 + \beta_0^\top X + \sum_{i,j=1}^p \beta_1^{(ij)} X^{(i)} X^{(j)} + \varepsilon.
\end{equation}
When applying our tests to this quadratic polynomial model, the resulting $p$-values of the tests $TCvM^{2}_{C}$, $TCvM^{2}_{CF}$, and $Hybrid_{CF}$ are $0.4759$, $0.6510$, and $0.5989$, respectively.
These results fail to reject the null hypothesis, implying that the polynomial regression model (\ref{quadratic-model}) may be plausible to fit this dataset. To further visualize this fit, we present a scatter plot of the residuals from the quadratic polynomial model versus the fitted values $\hat{Y}$ in Figure~\ref{fig:Res}. The absence of an obvious trend between the residuals and the fitted values further supports the adequacy of this model.

Next, we apply our proposed tests to the AML dataset, which was obtained from cBioPortal (\url{https://www.cbioportal.org/study/summary?id=aml_ohsu_2022}). This dataset consists of RNA-Seq expression profiles on $444$ patients, who are classified by the ELN2017 criteria into high-risk ($319$) and non–high-risk ($125$) groups. The expression data contains $22834$ genes, yielding a ultra-high dimensional setting for evaluation. We then evaluate the adequacy of a sparse linear logistic regression model for the AML classification task and apply our tests to check whether the functional form of the conditional expectation, $E[Y|X] = \frac{\exp(\beta^\top X)}{1 + \exp(\beta^\top X)}$, is plausible. The $p$-values of the tests $TCvM^{2}_{C}$, $TCvM^{2}_{CF}$, and $Hybrid_{CF}$ are approximately $0.9877$, $0.9961$, and $0.9921$, respectively. The universally high $p$-values lead us to fail to reject the null hypothesis, suggesting that the sparse linear logistic regression model is adequate for this dataset.
We further calculate the model predictive accuracy of the linear logistic regression model using $20$ repetitions of $5$-fold cross-validation. The resulting average predictive accuracy and Area Under the ROC Curve (AUROC) are $0.8514$ and $0.9292$, respectively. These metrics provide additional confirmation of the adequacy of the linear logistic regression model for this dataset.

\section{Discussion}
In this paper, we propose a two-step methodology for testing the goodness-of-fit of sparse parametric regression models, when the covariate dimension $p$ may significantly exceed the sample size $n$. In the first step, we construct the Cram\'{e}r-von Mises type test based on the martingale transformation of projected residual marked empirical processes. Under the null hypothesis, our projected tests are asymptotically distribution-free. Under the alternative hypothesis and mild conditions, the projected tests are consistent with asymptotic power $1$ for almost all projections on the unit sphere and can detect local alternatives departing from the null at the parametric rate of $O(n^{-1/2})$. In the second step, we employ the Cauchy combination method and data splitting to combine the projected tests to form our final tests, thereby enhancing power. Moreover, since empirical process-based tests are generally more sensitive to low-frequency alternatives and local smoothing tests are more powerful for high-frequency alternatives, we further propose a hybrid test that combines our empirical process-based tests and the local smoothing test proposed by \cite{tan2025b}. Simulation results show that the hybrid test performs very well for both low-frequency and high-frequency alternative models.
It is important to note that our methodology requires data splitting for the construction of the test statistics, which introduces variability in the values of the test statistics. An interesting question is whether the data splitting strategy can be completely avoided. Such alternative methods, without data splitting, would be particularly useful for model checking where there exist dependency between observations. We also note that model checking for the conditional mean function is a special case of testing conditional moment restrictions. It is of interest to extend our method to test general conditional moment restrictions in ultra-high dimensional settings.

%\renewcommand\refname{References}
%\nocite{*}
%\bibliographystyle{plain}
\bibliographystyle{Chicago}
\bibliography{mybib}

\end{document}